\documentclass[11pt]{amsart}
\usepackage{amssymb}
\usepackage[cmtip,matrix,arrow]{xy}

\newtheorem{theorem}{Theorem}[section]
\newtheorem{corollary}[theorem]{Corollary}
\newtheorem{proposition}[theorem]{Proposition}

\newtheorem{lemma}[theorem]{Lemma}

\theoremstyle{remark}
\newtheorem{remark}[theorem]{Remark}
\newtheorem{definition}[theorem]{Definition}
\newtheorem{remarks}[theorem]{Remarks}
\newtheorem{example}[theorem]{Example}
\newtheorem{examples}[theorem]{Examples}
\newcommand\A{\mathcal{A}}
\newcommand\be{\begin{equation}\label}
\newcommand\ee{\end{equation}}

\newcommand{\W}{\mathcal{W}}

\newcommand{\U}{\on{U}}

\newcommand{\E}{\mathcal{E}}

\newcommand{\R}{\mathbb{R}}
\newcommand{\C}{\mathbb{C}}
\newcommand{\F}{\mathbb{F}}

\newcommand{\Z}{\mathbb{Z}}

\newcommand{\Cl}{{\on{Cl}}}

\newcommand\lie[1]{\mathfrak{#1}}
\renewcommand{\k}{\lie{k}}
\newcommand{\h}{\lie{h}}
\renewcommand{\a}{\lie{a}}

\newcommand{\g}{\lie{g}}

\newcommand{\n}{\lie{n}}

\renewcommand{\t}{\lie{t}}

\newcommand{\so}{\lie{so}}
\newcommand{\on}{\operatorname}

\newcommand{\Ad}{ \on{Ad} } 
\newcommand{\ad}{ \on{ad} }

\newcommand{\End}{ \on{End} } 

\newcommand{\Hom}{ \on{Hom}}

\newcommand{\SO}{ \on{SO}} 
\renewcommand{\S}{\mathcal{S}}

\newcommand{\D}{ \mathcal{D} }

\newcommand\dirac{/\kern-1.2ex\partial} 
\newcommand\qu{/\kern-.7ex/} 

\renewcommand\a{\mathfrak{a}}
\newcommand{\lra}{\longrightarrow}
\newcommand{\hra}{\hookrightarrow}

\renewcommand{\d}{{\mbox{d}}}
\newcommand{\ol}{\overline}

\newcommand\sig{\sigma}
\newcommand\eps{\epsilon}
\newcommand\Om{\Omega}
\newcommand\om{\omega}

\newcommand{\f}{\frac}

\newcommand{\p}{\mf{p}}
\renewcommand{\l}{\langle}
\renewcommand{\r}{\rangle}
\newcommand{\hh}{{\textstyle \f{1}{2}}}
\newcommand{\ti}{\tilde}

\newcommand\beqn{\begin{equation}}      
\newcommand\eeqn{\end{equation}}      
\newcommand{\ca}{\mathcal}
\newcommand{\wh}{\widehat}
\newcommand{\wt}{\widetilde}
\newcommand{\mf}{\mathfrak}
\newcommand{\beq}{\begin{eqnarray*}}
\newcommand{\eeq}{\end{eqnarray*}}

\newcommand{\0}{{\bar{0}}}
\renewcommand{\1}{{\bar{1}}}

\newcommand{\da}{\on{da}}
\newcommand{\ds}{\on{ds}}
\newcommand{\dl}{\on{dl}}
\newcommand{\gds}{\g-\on{ds}}
\newcommand{\gda}{\g-\on{da}}
\newcommand{\gdl}{\g-\on{dl}}
\newcommand{\kds}{\k-\on{ds}}
\newcommand{\kda}{\k-\on{da}}
\newcommand{\kdl}{\k-\on{dl}}
\begin{document}

\title[Lie theory and the Chern-Weil homomorphism]
{Lie theory and the Chern-Weil homomorphism}

\author{A. Alekseev}
\address{University of Geneva, Section of Mathematics,
2-4 rue du Li\`evre, 1211 Gen\`eve 24, Switzerland}
\email{alekseev@math.unige.ch}

\author{E. Meinrenken}
\address{University of Toronto, Department of Mathematics,
100 St George Street, Toronto, Ontario M5S3G3, Canada }
\email{mein@math.toronto.edu}

\date{August 1, 2003}

\begin{abstract}
{Let $P \to B$ be a principal $G$-bundle. For any connection $\theta$
on $P$, the Chern-Weil construction of characteristic classes defines
an algebra homomorphism from the Weil algebra $W\g=S\g^*\otimes
\wedge\g^*$ into the algebra of differential forms 
$\A=\Om(P)$. Invariant
polynomials $(S\g^*)_{\on{inv}}\subset W\g$ map to cocycles, and the
induced map in cohomology $(S\g^*)_{\on{inv}}\to H(\A_{\on{basic}})$
is independent of the choice of $\theta$. The algebra $\Om(P)$ is
an example of a {\em commutative} $\g$-differential algebra with
connection, as introduced by H. Cartan in 1950. As
observed by Cartan, the Chern-Weil construction generalizes to all
such algebras.

In this paper, we introduce a canonical Chern-Weil map $W\g\to \A$ for
possibly {\em non-commutative} $\g$-differential algebras with connection. 
Our main observation is that the generalized Chern-Weil map
is an algebra homomorphism ``up to $\g$-homotopy''. 
Hence, the induced map $(S\g^*)_{\on{inv}} \to H_{\on{basic}}(\A)$ is an
algebra homomorphism. As in the standard Chern-Weil theory, 
this map is independent of the choice of connection. 

Applications of our results include: a conceptually
easy proof of the Duflo theorem for quadratic Lie algebras, a short 
proof of a conjecture of Vogan on Dirac cohomology,  generalized 
Harish-Chandra projections for quadratic Lie algebras, an extension 
of Rouvi\`{e}re's theorem for symmetric pairs, 
and a new construction of universal 
characteristic forms in the Bott-Shulman complex. }
\end{abstract}

\maketitle

\tableofcontents
\newpage

\section{Introduction}\label{sec:intro}
In an influential paper from 1950, H. Cartan \cite{ca:no} presented an
algebraic framework for the Chern-Weil \cite{ch:cha,we:oe1}
construction of characteristic classes in terms of differential forms.
In Cartan's approach, the de Rham complex $\Om(P)$ of differential
forms on a principal $G$-bundle $P$ is generalized to a differential
algebra $\A$, together with algebraic counterparts of the 
Lie derivative and contraction operations for the action
of the Lie algebra $\g$ of $G$. We
will refer to any such $\A$ as a {\em $\g$-differential
algebra}. Cartan introduced the notion of an {\em algebraic
connection} on $\A$; $\g$-differential algebras admitting connections
are called {\em locally free} and are viewed as algebraic counterparts
of principal bundles. A counterpart of the base of the
principal bundle is the basic subcomplex $\A_{\on{basic}}$.  The {\em
Weil algebra} $W\g=S\g^*\otimes\wedge\g^*$ replaces
the classifying bundle $EG\to BG$. The generators of
$\wedge\g^*$ are viewed as ``universal connections'', the generators
of $S\g^*$ as ``universal curvatures''. Cartan shows that if $\A$ is
any (graded) {\em commutative} $\g$-differential algebra with connection
$\theta$, there is a characteristic homomorphism
\begin{equation}\label{eq:char0}
 c^\theta:\,W\g\to \A
\end{equation}
sending the generators of $\wedge\g^*$ to the connection variables of
$\A$ and the generators of $S\g^*$ to the curvature variables of
$\A$. Passing to the cohomology of the basic subcomplex, this gives a
homomorphism
\begin{equation}\label{eq:chernweil} 
(S\g^*)_{\on{inv}}=H((W\g)_{\on{basic}})\to H(\A_{\on{basic}})
\end{equation}
from the algebra of invariant polynomials on the Lie algebra $\g$ into
the cohomology algebra of the basic subcomplex of $\A$. As in the usual
Chern-Weil theory, this homomorphism is independent of the choice of 
$\theta$.

The main theme of this paper is a generalization of Cartan's algebraic
Chern-Weil construction to possibly non-commu\-ta\-tive
$\g$-differential algebras. The idea in the general case is to define
$c^\theta$ by a suitable ordering prescription. Recall that any linear
map $E\to \A$ from a vector space to an associative algebra $\A$
extends to a linear map $S(E)\to \A$ from the symmetric algebra, by
symmetrization. This also holds for $\Z_2$-graded vector spaces and
algebras, using super symmetrization (i.e. taking signs into
account). The Weil algebra may be viewed as the (super) symmetric
algebra $S(E)$ over the space $E$ spanned by the generators
$\mu\in\g^*$ of $\wedge\g^*$ and their differentials $\ol{\mu}=\d\mu$,
and a connection $\theta$ on a $\g$-differential algebra $\A$ defines
a linear map $E\to \A$. With $c^\theta: W\g \to \A$ defined by
symmetrization, we prove:
\vskip.2in

\noindent{\bf Theorem A.}  {\it For any $\g$-differential algebra $\A$
with connection $\theta$, the map $c^\theta:\,W\g\to \A$ is a
homomorphism of $\g$-differential spaces. 
The induced homomorphism in
basic cohomology $(S\g^*)_{\on{inv}}\to H(\A_{\on{basic}})$ is
independent of $\theta$, and is an algebra homomorphism.}
\vskip.2in

In fact, we will show that any two homomorphisms of $\g$-differential spaces
$W\g\to \A$ induce the same map in basic cohomology, provided they 
agree on the unit element of $W\g$. Note that the Weil algebra 
could also be viewed as a super symmetric algebra over the subspace 
$E'\subset W\g$ spanned by the universal connections and curvatures.
However, the resulting symmetrization map $W\g=S(E')\to \A$ would not 
be a chain map, in general.

Our first application of Theorem A gives a new perspective on the
proof of \cite{al:no} of the Duflo isomorphism for quadratic Lie
algebras. Recall that a Lie algebra $\g$ is called {\em quadratic} if it
comes equipped with an invariant scalar product $B$. Let $\W\g$ be the
non-commutative super algebra generated by odd elements $\xi$ and even
elements $\ol{\xi}$ for $\xi\in\g$, subject to relations
$$ \bar{\xi}\bar{\xi}'-\bar{\xi}'\bar{\xi}=\ol{[\xi,\xi']_\g},\ \ 
\bar{\xi}\xi'-\xi'\bar{\xi}=[\xi,\xi']_\g,\ \
\xi\xi'+\xi'\xi=B(\xi,\xi').$$
Using $B$ to identify $\g^*$ and $\g$, we obtain a symmetrization map 
$W\g\to \W\g$. The following is a fairly easy consequence of Theorem A:
\vskip.2in 

\noindent {\bf Theorem B.} {\it There is a commutative diagram} 
$$ \xymatrix{ W\g\ar[r]& \W\g\\S\g\ar[r]\ar[u] & U\g\ar[u]} $$
{\it in which the vertical maps are injective algebra 
homomorphisms and the horizontal maps are vector space isomorphisms. 
The lower map restricts an algebra isomorphism on invariants.} 
\vskip.2in 

Recall that the Poincar\'{e}-Birkhoff-Witt symmetrization $S\g\to
U\g$ does {\em not} restrict to an algebra homomorphism on invariants, 
in general. On the other hand, it was shown by Duflo that the 
PBW map does have this property if it is pre-composed with a 
certain infinite order differential operator known as the 
``Duflo factor''.  
\vskip.2in 

\noindent {\bf Theorem C.}  {\it The lower horizontal map in the
commutative diagram of Theorem B coincides with the Duflo map
\cite{du:op}.}
\vskip.2in 

That is, while the Duflo map is not a symmetrization map for
$U\g$, it may be viewed as the {\em restriction} of a symmetrization
map of a super algebra containing $U\g$! We stress that our theory only
covers the case of quadratic Lie algebras -- it remains a mystery how
the general situation might fit into this picture.

Suppose $\k\subset\g$ is a quadratic subalgebra of $\g$, i.e. that 
the restriction of $B|_\k$ is non-degenerate. Let $\mf{p}$ denote the 
orthogonal complement to $\k$ in $\mf{p}$, and $\on{Cl}(\mf{p})$ 
its Clifford algebra. In \cite{ko:cu}, Kostant introduced a canonical 
element $\ca{D}_{\g,\k}$ of the algebra $(U\g\otimes 
\on{Cl}(\mf{p}))_{\k-\on{inv}}$ which he called the {\em cubic Dirac 
operator}. He showed that $\ca{D}_{\g,\k}$ squares to an element 
of the center of this algebra, so that the graded commutator 
$[\ca{D}_{\g,\k},\cdot]$ is a differential. The cohomology of 
this differential features in a conjecture of Vogan. Generalizing 
results of Huang-Pandzic \cite{hu:vo} and Kostant \cite{ko:di} 
we prove: 
\vskip.2in 

\noindent {\bf Theorem D.}  {\it There is a natural algebra
homomorphism $ (U\k)_{\k-\on{inv}}\to (U\g\otimes
\on{Cl}(\mf{p}))_{\k-\on{inv}},$ taking values in cocycles and
inducing an isomorphism in cohomology. The map 
$(U\g)_{\g-\on{inv}}\to (U\k)_{\k-\on{inv}}$ taking $z\in (U\g)_{\g-\on{inv}}$ 
to the cohomology class of $z\otimes 1$ coincides with the restriction map 
$(S\g)_{\g-\on{inv}}\to (S\k)_{\k-\on{inv}}$ under the Duflo isomorphisms
for $\g$ and $\k$. } 
\vskip.2in

Our next result is a Harish-Chandra map for a quadratic Lie algebra $\g$ 
with a decomposition $\g = \n_- \oplus \mf{k}\oplus \n_+$, where
$\k$ is a quadratic Lie subalgebra of $\g$ and $\n_\pm$ are 
$\k$-invariant isotropic subalgebras (that is, the restriction of
$B$ to $\n_\pm$ vanishes).
By the Poincar\'{e}-Birkhoff-Witt theorem, the
splitting of $\g$ gives rise to a
decomposition of the enveloping algebra $U\g$, hence to a projection
$\kappa_U:\,U\g\to U\k$.  As for the usual Harish-Chandra projection,
it is convenient to compose $\kappa_U$ with a certain automorphism
$\tau$ of $U\k$ (the ``$\rho$-shift'').
\vskip.2in 

\noindent {\bf Theorem E.} 
{\it Under the Duflo isomorphisms for $\g,\k$, the 
composition $\tau\circ \kappa_U:\,(U\g)_{\g-\on{inv}}\to 
(U\k)_{\k-\on{inv}}$ coincides with the projection 
$(S\g)_{\g-\on{inv}}\to (S\k)_{\k-\on{inv}}$. In particular, it is an 
algebra homomorphism.} 
\vskip.2in 

We obtain Theorem E by studying the Harish-Chandra projection
$\kappa_\W:\,\W\g\to \W\k$, and comparing to the natural projection
$\kappa_W:\,W\g\to W\k$. It turns out that $\kappa_\W$ directly
restricts to $\tau\circ \kappa_U$: That is, the shift 
$\tau$ emerges from the theory in a very natural way and need not 
be put in 'by hand'. 

Let $(\g,\k)$ be a symmetric pair, that is, $\k$ is the fixed point
set of an involutive automorphism $\epsilon$ on $\g$.  Let $\mf{p}$ be
its complement given as the eigenspace of $\epsilon$ for the
eigenvalue $-1$.  By results of Lichnerowicz \cite{li:op1} and Duflo
\cite{du:ope}, the algebra $(U\g/U\g\k^f)_{\k-\on{inv}}$ (where
$\k^f\hra U\g$ a suitable ``twisted'' inclusion of $\k$) is {\em
commutative}. Rouvi\`{e}re in his paper \cite{ro:es} introduced a map
$(S\mf{p})_{\k-\on{inv}}\to (U\g/U\g\k^f)_{\k-\on{inv}}$ generalizing
the Duflo isomorphism, and described conditions under which this map 
is an algebra isomorphism. We prove a similar result for the 
following new class of examples:
\vskip.2in 

\noindent {\bf Theorem F.}  {\it Suppose $\g$ carries an invariant
scalar product $B$ that changes sign under $\epsilon$. Then the
Duflo-Rouvi\`{e}re map $(S\mf{p})_{\k-\on{inv}}\to
(U\g/U\g\k^f)_{\k-\on{inv}}$ is an algebra isomorphism. }
\vskip.2in 

Anti-invariance of $B$ under the involution $\eps$ implies that $B$
vanishes on both $\k$ and $\mf{p}$, and gives a non-degenerate pairing
between the two subspaces. In line with our general strategy, we prove
this result by identifying the Duflo-Rouvi\`{e}re map as a Chern-Weil
map, using the isomorphism
$(S\mf{p})_{\k-\on{inv}}=(S\k^*)_{\k-\on{inv}}$ given by the pairing.
We also describe a generalization to 
isotropic subalgebras $\k\subset\g$ of quadratic Lie algebras.

Our final result is a new construction of universal characteristic
forms in the Bott-Shulman complex. For any Lie group $G$,
Bott and Shulman considered a double complex $\Om^q(G^p)$ as a model for
differential forms on the classifying space $BG$, and showed how to
associate to any invariant polynomial on $\g$ a cocycle for the total
differential on this double complex.  In our alternative approach, we
observe that $\bigoplus_{p,q} \Om^q(G^p)$ carries a 
non-commutative product, and obtain: 
\vskip.2in 

\noindent {\bf Theorem G.} 
{\it The generalized Chern-Weil construction defines a linear map 
$$ (S\g^*)_{\on{inv}}\to \bigoplus_{p,q} \Om^q(G^{p})$$
taking values in cocycles for the total differential. The image of a
polynomial of degree $r$ under this map has non-vanishing components
only in bidegree $p+q=2r$ with $p\le r$. The map induces an algebra
homomorphism in cohomology, and in fact an algebra isomorphism if $G$
is compact.}
\vskip.1in

\noindent{\bf Acknowledgments.}  We are grateful to A. Beilinson,
V. Drinfeld, M. Duflo, E. Getzler, V. Ginzburg, B. Kostant, S. Kumar,
N. Markarian, J. Roberts, P. Severa, V. Turaev and M. Vergne for their
interest in our work and for useful comments. This work was supported
in part by NSERC, by the Swiss National Science Foundation and by the Erwin
Schr\"odinger Institute for Mathematical Physics.


\section{Non-commutative differential algebras}\label{sec:nonco}
In this Section we review some material on symmetrization maps for
super vector spaces and $\g$-differential spaces. Our conventions for super
spaces will follow \cite{del:no}; in particular we take the
categorical point of view that super vector spaces form a tensor category
where the super sign convention is built into the isomorphism
$E\otimes E'\to E'\otimes E$. The concept of $\g$-differential spaces
is due to Cartan \cite{ca:no}, a detailed treatment can be found in
the book \cite{gu:su}.

\subsection{Conventions and notation} 
Throughout, we will work over a
field $\F$ of characteristic $0$. A super vector space is a vector space over
$\F$ with a $\Z_2$-grading $E=E^\0\oplus E^\1$. Super vector spaces form an
$\F$-linear tensor category; algebra objects in this tensor category
are called super algebras, Lie algebra objects are called super Lie
algebras. If $E,E'$ are super vector spaces, we denote by $L(E,E')$ the super
space of all linear maps $A:\,E\to E'$, and by $\on{End}(E)=L(E,E)$
the super algebra of endomorphisms of $E$.  By contrast, the space of
homomorphisms $\on{Hom}(E,E')=L(E,E')^\0$ consists of only the {\em
even} linear maps.

\subsection{Symmetrization maps} \label{subsec:symm}
Let $E=E^\0\oplus E^\1$ be a super vector space.  The (super) symmetric
algebra $S(E)=\bigoplus_{k=0}^\infty S^k(E)$ is the quotient of the
tensor algebra $\ca{T}(E)=\bigoplus_{k=0}^\infty E^{\otimes{k}}$ by
the two-sided ideal generated by all elements of the form $v\otimes
w-(-1)^{|v||w|}w\otimes v$, for homogeneous elements $v,w\in E$ of
$\Z_2$-degree $|v|,|w|$. Both $\ca{T}(E)$ and $S(E)$ are super
algebras, in such a way that the inclusion of $E$ is a homomorphism of
super vector spaces. The tensor algebra $\ca{T}(E)$ is characterized by the
universal property that any homomorphism of super vector spaces $E\to \A$
into a super algebra $\A$ extends uniquely to a homomorphism of super
algebras $\ca{T}(E)\to \A$; the symmetric algebra has a similar
universal property for commutative super algebras.

Given a super algebra $\A$, any homomorphism of super vector spaces 
$\phi:\,E\to \A$ extends to $S(E)$ by symmetrization 
$$ \on{sym}(\phi):\,S(E)\to \A,\ \ 
v_1\cdots v_k\mapsto \f{1}{k!} \sum_{\sig\in \mathfrak{S}_k}
(-1)^{N_\sig(v_1,\ldots,v_k)} \phi(v_{\sig^{-1}(1)}) \cdots \phi(v_{\sig^{-1}(k)}).$$
Here $\mathfrak{S}_k$ is the symmetric group, and
${N_\sig(v_1,\ldots,v_k)}$ is the number of pairs $i<j$ such that
$v_i,v_j$ are odd elements and $\sig^{-1}(i)>\sig^{-1}(j)$.

Equivalently, the symmetrization map $S(E)\to \A$ may be characterized as 
the inclusion $S(E)\to \ca{T}(E)$ as ``symmetric tensors'', followed by 
the algebra homomorphism $\ca{T}(E)\to \A$ given by the universal property of 
$\ca{T}(E)$.

\subsection{Poincar\'{e}-Birkhoff-Witt symmetrization}
If $(E,[\cdot,\cdot]_E)$ is a super Lie algebra, one defines the
enveloping algebra $U(E)$ as the quotient of $\ca{T}(E)$ by the
relations $v_1\otimes v_2-(-1)^{|v_1||v_2|}v_2\otimes
v_1-[v_1,v_2]_E$.  By the Poincar\'{e}-Birkhoff-Witt theorem for super Lie
algebras (Corwin-Neeman-Sternberg \cite{co:gr}, see also
\cite{del:no}) the symmetrization map $S(E)\to U(E)$ is a linear
isomorphism.
 
Similarly, if $E$ is a super vector space with a skew-symmetric bi-linear form $\om\in
\Hom(E\otimes E,\F)$ 
(i.e. $\om(v,w)=-(-1)^{|v||w|}\om(w,v)$), one defines the {\em Weyl
algebra} $\on{Weyl}(E)$ as the quotient of the tensor algebra by the
ideal generated by elements $v\otimes w-(-1)^{|v||w|}w\otimes
v-\om(v,w)$. The corresponding symmetrization map $S(E)\to
\on{Weyl}(E)$ is an isomorphism of super vector spaces, known as 
the {\em quantization map} for the Weyl algebra. 
\footnote{The fact that $q$ is an isomorphism may be deduced from the
PBW isomorphism for the Heisenberg Lie algebra $E\oplus\F\mf{c}$,
i.e. the central extension with bracket $[v,v']=\om(v,v')\mf{c}$.
Indeed, the symmetrization map $S(E\oplus\F\mf{c})\to
U(E\oplus\F\mf{c})$ restricts to an isomorphism between the ideals
generated by $\mf{c}-1$, and so the claim follows by taking quotients
by these ideals.} In the purely odd case $E^\0=0$, $\om$ is a
symmetric bilinear form $B$ on $V=E^\1$ (viewed as an ungraded vector
space), the Weyl algebra is the Clifford algebra of $(V,B)$, and the
symmetrization map reduces to the Chevalley quantization map
$q:\,\wedge V\to \on{Cl}(V)$.

\subsection{Derivations}\label{subsec:der}
Given a super algebra $\A$ we denote by $\on{Der}(\A)\subset \End(\A)$
the super Lie algebra of derivations of $\A$. Similarly, if $E$ is a
super Lie algebra we denote by $\on{Der}(E)\subset
\End(E)$ the super Lie algebra of derivations of $E$.
For any super vector space $E$ there is a unique homomorphism of super
Lie algebras
$$\on{End}(E) \to\on{Der}(\ca{T}(E)),\,A\mapsto D_A^{\ca{T}(E)}$$
such that $D_A^{\ca{T}(E)}(v)=Av$ for $v\in E\subset \ca{T}(E)$.
Similarly one defines $D_A^{S(E)}\in \on{Der}(S(E))$ and, if $A$ is a
derivation for a super Lie bracket on $E$, $D_A^{U(E)}\in \on{Der}(U(E))$.
We will need the following elementary fact. 
\begin{lemma}\label{lem:deriv}
Let $E$ be a super vector space, $\A$ a super algebra, and $\phi\in
\on{Hom}(E,\A)$ a homomorphism of super vector spaces. Suppose we are
given a linear map $A\in \End(E)$ and a derivation $D\in
\on{Der}(\A)$, and that $\phi$ intertwines $A$ and $D$. Then the
extended map $\on{sym}(\phi): \,S(E)\to \A$ obtained by symmetrization
intertwines $D_A^{S(E)} \in \on{Der}(S(E))$ and $D$.
\end{lemma}
\begin{proof}
Recall that $\on{sym}(\phi)$ factors through the symmetrization map 
for the tensor algebra $\ca{T}(E)$. Since the map $\ca{T}(E)\to \A$
intertwines $D_A^{\ca{T}(E)}$ with $D$, it suffices to prove
the Lemma for $\A=\ca{T}(E)$, $D=D^{\ca{T}(E)}$.  The action of
$D_A^{\ca{T}(E)}$ on $E^{\otimes k}$ commutes with the action of
$\mathfrak{S}_k$, as one easily checks for transpositions. In
particular $D_A^{\ca{T}(E)}$ preserves the $\mathfrak{S}_k$-invariant
subspace. It therefore restricts to $D_A^{S(E)}$ on $S(E)\subset
\ca{T}(E)$.
\end{proof}

\subsection{Differential algebras}\label{subsec:da}
A {\em differential space} ($\ds$) is a super vector space $E$, together with
a differential, i.e. an odd endomorphism $\d\in \on{End}(E)^\1$
satisfying $\d\circ\d=0$. Morphisms in the category of differential
spaces will be called chain maps or {\em $\ds$ homomorphisms}. The
tensor product $E\otimes E'$ of two differential spaces is a
differential space, with $\d(v\otimes v')=\d v\otimes
v'+(-1)^{|v|}v\otimes \d v'$. Algebra objects in this tensor category
are called {\em differential algebras} ($\da$),
Lie algebra objects are called {\em differential Lie algebras}
($\dl$).  The discussion from Section \ref{subsec:der} shows:

\begin{lemma}
Let $\A$ be a differential algebra, and $\phi:\,E\to \A$ 
a $\ds$ homomorphism. Then the symmetrized map $S(E)\to \A$ 
is a $\ds$ homomorphism. 
\end{lemma}

For any differential algebra $(\A,\d)$, the unit $i:\,\F\hra \A$
(i.e. the inclusion of $\F$ as multiples of the unit element) is a
$\da$ homomorphism. By an {\em augmentation map} for $(\A,\d)$, we
mean a $\da$ homomorphism $\pi:\,\A\to \F$ such that $\pi\circ i$ is
the identity map.  For the tensor algebra, the natural projections
onto $\F=E^{\otimes 0}$ is an augmentation map; this descends to
augmentation maps for $S(E)$ and (in the case of a super Lie algebra)
$U(E)$.

\subsection{Koszul algebra}\label{subsec:koszul}
For any vector space $V$ over $\F$, let $E_V$ be the $\ds$  
with $E_V^\0=V$ and $E_V^\1=V$, with differential $\d$ 
equal to $0$ on even elements and given by the natural isomorphism 
$E_V^\1\to E_V^\0$ on odd elements. For $v\in V$ denote the corresponding 
even and odd elements in $E_V$ by $\ol{v}\in E_V^\0$ and 
$v\in E_V^\1$, respectively. Thus 
$$ \d v=\ol{v},\ \ \d\ol{v}=0.$$
The symmetric algebra $S(E_V)$ is known as the {\em Koszul algebra}
over $V$. It is characterized by the universal property that if $\A$
is any commutative $\da$, any vector space homomorphism $V\to \A^\1$
extends to a unique homomorphism of $\da$'s $S(E_V)\to \A$.  We will
also encounter a non-commutative version of the Koszul algebra,
$\ca{T}(E_V)$. It has a similar universal property, but in the
category of not necessarily commutative $\da$'s $\A$.  The
non-commutative Koszul algebra appears in a
paper of Gelfand-Smirnov \cite{ge:al}. 

For any Lie algebra $\g$, we denote by $\wt{\g}$ the super Lie algebra
$$ \wt{\g}=\g\ltimes\g$$
(semi-direct product) where the even part $\wt{\g}^\0=\g$ acts on the
odd part $\wt{\g}^\1=\g$ by the adjoint representation. It is a differential 
Lie algebra under the identification $\wt{\g}=E_\g$: that is,  
$\d$ is a derivation for the Lie bracket. 

\subsection{$\g$-differential spaces}
A {\em $\g$-differential space} ($\gds$) is a differential space $(E,\d)$
together with a $\dl$ homomorphism, $\wt{\g}\to \on{End}(E)$. 
That is, it consists of a representation of the super Lie algebra 
$\wt{\g}$ on $E$, where the operators $\iota_\xi,L_\xi\in 
\on{End}(E)$ corresponding to $\xi,\ol{\xi}\in \wt{\g}$ satisfy the 
relations
\begin{equation}\label{eq:defrelations}
 [\d,\iota_\xi]=L_\xi,\ \ [\d,L_\xi]=0.
\end{equation}
The operators $\iota_\xi$ are called {\em contractions}, the operators
$L_\xi$ are called {\em Lie derivatives}.  The tensor product of any
two $\gds$, taking the tensor product of the $\wt{\g}$-representations,
is again a $\gds$. Hence $\gds$'s form an $\F$-linear tensor category;
the algebra objects in this tensor category are called {\em
$\g$-differential algebras} ($\gda$), the Lie algebra objects are
called {\em $\g$-differential Lie algebras} ($\gdl$). This simply
means that the representation should act by derivations of the
product respectively Lie bracket.

For any $\gds$ $E$, one defines the horizontal subspace
$E_{\on{hor}}=\bigcap \on{ker}\iota_\xi$, the invariant subspace
$E_{\on{inv}}=\bigcap \on{ker}L_\xi$, and the basic subspace
$E_{\on{basic}}=E_{\on{hor}}\cap E_{\on{inv}}$. That is,
$E_{\on{basic}}$ is the space of fixed vectors for
$\wt{\g}$. This subspace is stable under $\d$, hence is a
differential space.  Any $\gds$ homomorphism $\phi:\,E\to E'$
restricts to a chain map between basic subspaces.

An example of a $\gdl$ is $E=\wt{\g}$ with the adjoint action. 
Another example is obtained by adjoining $\d$ as an odd element, 
defining a semi-direct product
\begin{equation}  \F\d\ltimes \wt{\g}  \label{eq:hatg} \end{equation}
where the action of $\d$ on $\wt{\g}$ is as the differential, $\d\xi=\ol{\xi},\ 
\d\ol{\xi}=0$. (Note that a $\gds$ can be
defined equivalently as a module for the super Lie algebra
\eqref{eq:hatg}; this is the point of view taken in the book
\cite{gu:su}.) The symmetric and tensor algebras over a $\gds$ $E$
are $\gda$'s, and Lemma \ref{lem:deriv} shows: 
\begin{lemma}\label{lem:gds}
If $\A$ is a $\gda$, and $\phi:\,E\to \A$ is a $\gds$ homomorphism, 
then the symmetrization $\on{sym}(\phi):\,S(E)\to \A$ is 
a $\gds$  homomorphism.
\end{lemma}

\subsection{Homotopy operators} \label{subsec:ho}
The space $L(E,E')$ of linear maps $\phi:\,E\to E'$ between 
differential spaces is itself a differential space, with differential 
$\d(\phi)=\d\circ \phi-(-1)^{|\phi|}\phi\circ \d$. 
Chain maps correspond to cocycles in $L(E,E')^\0$. 

A homotopy operator between two chain maps $\phi_0,\,\phi_1:\,E\to E'$ 
is an odd linear map $h\in L(E,E')^\1$ such that
$\d(h)=\phi_0-\phi_1$. A homotopy inverse to a chain map $\phi:\,E\to
E'$ is a chain map $\psi:\,E'\to E$ such that $\phi\circ\psi$ and
$\psi\circ\phi$ are homotopic to the identity maps of $E',E$.

\begin{lemma}\label{lem:homotopy}
Let $(E,\d)$ be a differential space, and suppose there exists
$s\in\on{End}(E)^\1$ with $[\d,s]=\on{id}_E$. Then the inclusion 
$i:\F\to \ca{T}(E)$ and the augmentation map $\pi:\,\ca{T}(E)\to 
\F$ are homotopy inverses. Similarly for the symmetric algebra. 
\end{lemma}

\begin{proof}
The derivation extension of $[\d,s]$ to $\ca{T}(E)$ is the Euler
operator on $\ca{T}(E)$, equal to $k$ on $E^{\otimes k}$.  Hence 
$[\d,s]+i\circ \pi\in \on{End}(E)^\0$ is an invertible chain map, and 
the calculation
$$ I-i\circ \pi=[\d,s]([\d,s]+i\circ \pi)^{-1}
=[\d,s([\d,s]+i\circ \pi)^{-1}]
$$ 
shows that $h=s([\d,s]+i\circ \pi )^{-1}\in \on{End}(E)^\1$
is a homotopy operator between $I$ and $i\circ \pi$. 
The proof for the symmetric algebra is similar.
\end{proof}

The Lemma shows that the Koszul algebra $S(E_V)$ and its non-commutative
version $\ca{T}(E_V)$ have trivial cohomology, since
$s\in\on{End}(E_V)^\1$ given by $s(\ol{v})=v$ and $s(v)=0$ has the
desired properties. We will refer to the corresponding $h$ as the
standard homotopy operator for the Koszul algebra.

More generally, if $E,E'$ are $\gds$ then the space $L(E,E')$ inherits
the structure of a $\gds$, with contractions
$\iota_\xi(\phi)=\iota_\xi\circ \phi- (-1)^{|\phi|}\phi\circ
\iota_\xi$ and Lie derivatives $L_\xi(\phi)=L_\xi\circ \phi-\phi\circ
L_\xi$.  $\gds$ homomorphisms $E\to E'$ are exactly the cocycles in
$L(E,E')_{\on{basic}}^\0$. A homotopy between two $\gds$ homomorphisms
$\phi_0,\phi_1:\,E\to E'$ is called a {\em $\g$-homotopy} if it is
$\wt{\g}$-equivariant.  Note that any such $h$ restricts to a homotopy of
$\phi_0,\phi_1:\,E_{\on{basic}}\to E'_{\on{basic}}$.

\subsection{Connection and curvature}
A {\em connection} on a $\gda$ $\A$ is a linear map
$\theta:\,\g^*\to \A^\1$ with the properties,
$$ \iota_\xi(\theta(\mu))=\mu(\xi),\ \ L_\xi(\theta(\mu))=
-\theta(\ad_\xi^*\mu).$$
A $\gda$ admitting a connection is called {\em locally free}. 
The {\em curvature} of a connection $\theta$ is the linear map
$F^\theta:\,\g^*\to \A^\0$ defined by
$$ F^\theta=\d\theta+\hh [\theta,\theta].$$ 
Here $[\theta,\theta]$ denotes the composition $ \g^*\to
\g^*\otimes\g^*\stackrel{\theta\otimes\theta}{\lra} \A^\1\otimes\A^\1\to
\A^\0,$ where the first map is the map dual to the Lie bracket and the
last map is algebra multiplication. The curvature map is equivariant
and satisfies $\iota_\xi F^\theta=0$.


The following equivalent definition of a connection will be useful in 
what follows. Let $\F\mf{c}$ be the 1-dimensional space spanned by an
even generator $\mf{c}$, viewed as a differential space on which $\d$
acts trivially, and consider the direct sum $E_{\g^*}\oplus \F\mf{c}$
with $\wt{\g}$-action given by 
\begin{equation}\label{eq:gdstructure}
 L_\xi\ol{\mu}=-\ol{\ad_\xi^*\mu},\ 
L_\xi\mu=-\ad_\xi^*\mu,\ 
\iota_\xi\ol{\mu}=-\ad_\xi^*\mu,\   
\iota_\xi\mu=\mu(\xi)\mf{c}.
\end{equation}
Then, a connection on a $\gda$ $\A$ is equivalent to a 
$\gds$ homomorphism 
\begin{equation}\label{eq:conn}
E_{\g^*}\oplus \F\mf{c}\to \A
\end{equation}
taking $\mf{c}$ to the unit of $\A$. 
\begin{remark}
It may be verified that $E_{\g^*}\oplus \F\mf{c}$ is the {\em odd} 
dual space of the super Lie algebra $\F\d\ltimes \wt{\g}$, 
i.e. the dual 
space with the opposite $\Z_2$-grading.   
\end{remark}

\section{The Chern-Weil homomorphism}\label{sec:weil}
\subsection{The Weil algebra}
The {\em Weil algebra} $W\g$ is a commutative $\gda$ with connection
$\g^*\to W\g$, with the following universal property: For any
commutative $\gda$ $\A$ with connection $\theta:\,\g^*\to \A$, there
exists a unique $\gda$ homomorphism $c^\theta:\, W\g\to \A$ such
that the following diagram commutes:
$$
\xymatrix{W\g\ar[r]_{c^\theta}&\A\\
\g^*\ar[u]\ar[ru]_{\theta}& }$$
We will refer to $c^\theta$ as the {\em characteristic homomorphism}
for the connection $\theta$.  The Weil algebra is explicitly given as 
a quotient 
\begin{equation}\label{eq:weil}
 W\g=S(E_{\g^*}\oplus \F\mf{c})/<\mf{c}-1>,\end{equation}
where $<\mf{c}-1>$ denotes the two-sided ideal generated by
$\mf{c}-1$. From the description \eqref{eq:conn} of connections, it is
obvious that $W\g$ carries a ``tautological'' connection. If $\A$ is
commutative, the homomorphism \eqref{eq:conn} extends, by the
universal property of the symmetric algebra, to a $\gda$ homomorphism
\begin{equation}\label{eq:conn1}
S(E_{\g^*}\oplus \F\mf{c}) \to \A
\end{equation}
This homomorphism takes $<\mf{c}-1>$ to $0$, and therefore descends to
a $\gda$ homomorphism $c^\theta:\,W\g\to \A$. If $\A$ is a {\em
non-commutative} $\gda$ with connection $\theta$, the map \eqref{eq:conn1} is
still well-defined, using symmetrization. Hence one obtains a 
canonical $\gds$ homomorphism 
\begin{equation}\label{eq:ch}
 c^\theta:\ W\g\to \A
\end{equation}
even in the non-commutative case. Of course, \eqref{eq:ch} is no longer an 
algebra homomorphism in general. What we will show below is that 
\eqref{eq:ch} is an algebra homomorphism ``up to $\g$-homotopy''. 

As a differential algebra, the Weil algebra is just the Koszul algebra
$W\g=S(E_{\g^*})$. The $\gda$-structure is given on generators
$\mu,\ol{\mu}$ by formulas similar to \eqref{eq:gdstructure}, with
$\mf{c}$ replaced by $1$, and the connection is the map sending
$\mu\in\g^*$ to the corresponding odd generator of $W\g$. 
The curvature map for this connection is given by 
$$  F:\,\g^*\to W\g,\ \mu\mapsto \wh{\mu}:=\ol{\mu}-\lambda(\mu)$$
where 
$$  \lambda:\,\g^*\to \wedge^2\g^*,\ 
\iota_\xi\lambda(\mu)=-\ad_\xi^*\mu,\ \ \xi\in\g,\ \mu\in\g^*$$
is the map dual to the Lie bracket, and we identify $\wedge\g^*$ as the subalgebra of $W\g$ defined by the odd generators.
Recall $\iota_\xi\wh{\mu}=0$ since the curvature of a connection is
horizontal.  The curvature map extends to an algebra homomorphism
$S\g^*\to (W\g)_{\on{hor}}$, which is easily seen to be an
isomorphism. Thus
\begin{equation}\label{eq:factWg}
 W\g=S\g^*\otimes\wedge\g^*,
\end{equation}
where $S\g^*\cong (W\g)_{\on{hor}}$ is generated by the ``curvature
variables '' $\wh{\mu}$, and the exterior algebra $\wedge\g^*$ by the
``connection variables'' $\mu$.  The Weil
differential vanishes on $(W\g)_{\on{basic}}=(S\g^*)_{\on{inv}}$;
hence the cohomology of the basic subcomplex coincides with
$(S\g^*)_{\on{inv}}$.

\begin{remark}\label{rem:wg1}
Let $\ti{W}\g$ be the non-commutative $\gda$, defined similar to
\eqref{eq:weil} but with the tensor algebra in place of the symmetric
algebra,
$$\ti{W}\g=\ca{T}(E_{\g^*}\oplus \F\mf{c})/<\mf{c}-1>.$$
Then $\ti{W}\g$ has a universal property similar to $W\g$ among 
{\em non-commutative} $\gda$'s with connection.
As a differential algebra, $\ti{W}\g$ is the non-commutative Koszul
algebra, $\ti{W}\g=\ca{T}(E_{\g^*})$. 
In particular, it is acyclic, with a canonical homotopy operator.
\end{remark}

\subsection{The rigidity property of the Weil algebra}
The algebraic version of the Chern-Weil theorem asserts that if $\A$
is a commutative locally free $\gda$, the characteristic homomorphisms
corresponding to any two connections on $\A$ are $\g$-homotopic (see
e.g. \cite{gu:su}).  This then implies that the induced map
$(S\g^*)_{\on{inv}}\to H(\A_{\on{basic}})$ is independent of the
choice of connection. The following result generalizes the Chern-Weil
theorem to the non-commutative setting:
\begin{theorem}[Rigidity Property]\label{th:rigidity}
Let $\A$ be a locally free $\gda$. 
Any two $\gds$ homomorphisms $c_0,c_1:\,W\g\to \A$ that agree 
on the unit of $W\g$ are $\g$-homotopic. 
\end{theorem}
The proof will be given in Section \ref{subsec:proof}. 
We stress that the maps $c_i$ need not be algebra homomorphism, 
and are not necessarily characteristic homomorphism for connections 
on $\A$. 
The rigidity theorem implies the following surprising result.
\begin{corollary}\label{cor:algebra}
Let $\A$ be a locally free $\gda$ and $c:\,W\g\to \A$ be a $\gds$ homomorphism
taking the unit of $W\g$ to the unit of $\A$. Then the induced map in basic cohomology
\begin{equation}\label{eq:alghom}
(S\g^*)_{\on{inv}}\to H(\A_{\on{basic}})\end{equation}
is an algebra
homomorphism, independent of the choice of $c$.
\end{corollary}
In particular, Corollary \ref{cor:algebra} applies to 
characteristic homomorphisms $c^\theta$ for connections $\theta$ 
on $\A$.
\begin{proof}
By Theorem \ref{th:rigidity}, the map 
\eqref{eq:alghom} is independent of $c$. Given
$w'\in (W\g)_{\on{basic}}=(S\g^*)_{\on{inv}}$, the two 
$\gds$ homomorphisms 
$$W\g\to \A,\,w\mapsto c(w w'),\ \ \  w\mapsto c(w)c(w').$$
agree on the unit of $W\g$, and are therefore   $\g$-homotopic. Hence
$[c(w w')]=[c(w) c(w')] = [c(w)] \,[c(w')]$ so that \eqref{eq:alghom} 
is an algebra homomorphism.
\end{proof}

\begin{remark}
Theorem \ref{th:rigidity} becomes false, in general, if one drops the
assumption that $\A$ is locally free.  To construct a counter-example,
let $\A_0$ and $\A_1$ be two $\Z$-graded locally free $\gda$'s with
characteristic maps $c_i$ for some choice of connections. Assume
that there is an element $p \in (S^k\g^*)_{\on{inv}}$ with $k \geq 2$
such that $c_0(p) \neq 0$ while $c_1(p)=0$.  (For instance, one can
choose $\A_0= W\g$, and $\A_1=\wedge \g^*$ the quotient of $W\g$ by
the ideal generated by the universal curvatures. Then the assumption
holds for any element $p$ of degree $k \geq 2$, taking $c_i$ to be the
characteristic maps for the canonical connections.) Let $E=(\A_0
\oplus \A_1)/V$, where $V$ is the one dimensional subspace spanned by
$(1_{\A_0} - 1_{\A_1})$.  View $E$ as an algebra with zero product,
and define a $\gda$ $\A= \F \oplus E$ by adjoining a unit.  The maps
$c_i$ descend to $\gds$ homomorphisms $\tilde{c}_i: W\g \to E\subset
\A$ that agree on $1_{W\g}$, while $\tilde{c}_0(p) \neq 0$ and
$\tilde{c}_1(p)=0$. Hence, $\tilde{c}_i$ are not $\g$-homotopic.
\end{remark}

\subsection{Proof of Theorem \ref{th:rigidity}}\label{subsec:proof}
The proof of Theorem \ref{th:rigidity} will be based on the 
following two observations. 
\vskip.1in
\begin{enumerate}
\item The Weil algebra $W\g$ fits into a family of $\g$-differential algebras
$S(E_{\g^*}\oplus \F\mf{c})/<\mf{c}-t>$ for $t\in\F$. As differential
algebras, these are all identified with $S(E_{\g^*})$.  The formulas
for the Lie derivatives $L_\xi$ are the same as for $W\g$, while the
contraction operators for the deformed algebras read
$$
\iota_\xi^t\ol{\mu}=-\ad_\xi^*\mu,\ \ 
\iota_\xi^t\mu=t\l\mu,\xi\r.
$$
It is clear that the resulting $\gda$'s for $t\not=0$ are all
isomorphic to $W\g$, by a simple rescaling. However, for $t=0$ one
obtains a non-isomorphic $\gda$ which may be identified with the
symmetric algebra over the space $E_{\g^*}$ with the co-adjoint
representation of $\wt{\g}$. The standard homotopy operator $h$ for 
$S(E_{\g^*})$ becomes a $\g$-homotopy operator
exactly for $t=0$.
\item  $S(E_{\g^*})$ is a
super Hopf algebra, with co-multiplication 
$\Delta:\,S(E_{\g^*})\to S(E_{\g^*})\otimes S(E_{\g^*})$ 
induced by the diagonal embedding $E_{\g^*}\to E_{\g^*}\oplus
E_{\g^*}$, and co-unit the augmentation map $\pi:\,S(E_{\g^*})\to \F$.
Both $\d$ and $L_\xi$ are co-derivations for this co-product.
For the contraction operators $\iota_\xi^t$ one finds, 
$$\Delta \circ \iota_\xi^{t_1+t_2}=
 (\iota_\xi^{t_1}\otimes 1+1\otimes \iota_\xi^{t_2})\circ \Delta ,
$$
in particular, $\iota_\xi^t$ is a co-derivation only for $t=0$.  
\end{enumerate}

\begin{proof}[Proof of Theorem \ref{th:rigidity}]
View the space $L(W\g,\A)$ of linear maps as a $\gds$, as 
in Section \ref{subsec:ho}. 
We have to show that if $c:\,W\g\to \A$ is 
a $\gds$ homomorphism with $c(1)=0$, then 
$c=\d(\psi)$ for a basic element $\psi\in L(W\g,\A)$.
Letting $h$ be the standard homotopy operator for $W\g=S(E_{\g^*})$, the 
composition $c\circ h$ satisfies $\d(c\circ h)=c$, however, it is not basic 
since $h$ is not a $\g$-homotopy.  

To get around this problem, we will exploit that 
the space $L(W\g,\A)$ has an algebra structure, using the
co-multiplication $\Delta$ on $W\g$ and the co-unit $\pi:\,W\g\to \F$. The
product of two maps $\phi_1,\phi_2:\,W\g\to \A$ is given by the
composition
$$ \phi_1\cdot\phi_2:\, W\g\stackrel{\Delta }{\lra} W\g\otimes
W\g\stackrel{\phi_1\otimes\phi_2}{\lra} \A\otimes \A\to \A,$$
and the unit element in $L(W\g,\A)$ is the composition 
$$1_L:=i_\A\circ \pi:\,W\g\to \A,$$
where $i_\A:\,\F\to \A$ is the unit for the algebra $\A$.  
The deformed contraction operators $\iota_\xi^t$ on
$W\g$ induce deformations of the contraction operators for
$L(W\g,\A)$, denoted by the same symbol. The formula in (b) implies
that
\begin{equation}\label{eq:prod12}
 \iota_\xi^{t_1+t_2}(\phi_1\cdot\phi_2)=
(\iota_\xi^{t_1}\phi_1) \cdot\phi_2+(-1)^{|\phi_1|}\phi_1\cdot 
\iota_\xi^{t_2}\phi_2.
\end{equation}
The sum $\phi=1_L+c$ is a $\gds$ homomorphism taking the unit in $W\g$
to the unit in $\A$. Since $c\circ i_{W\g}=0$, $\phi$ is an invertible
element of the algebra $L(W\g,\A)$, with inverse given as
a geometric series $\phi^{-1}=\sum_N (-1)^N c^N$.
The series is well-defined, since $c$ is locally
nilpotent: Indeed, $c^N$ vanishes on $S^j(E_{\g^*})$ for $N>j$.

Since $\iota_\xi^1(\phi)=0$ the inverse satisfies
$\iota_\xi^{-1}(\phi^{-1})=0$; this follows from \eqref{eq:prod12}
since the unit element $1_L \in L(W\g,\A)$ is annihilated by
$\iota_\xi^0$. Hence, the product $c \cdot \phi^{-1}\in L(W\g,\A)$ is
annihilated by $\d,L_\xi,\iota_\xi^0$. We now set, 
$$ \psi:= ((c \cdot \phi^{-1})\circ h)\cdot \phi.$$
Recall that $h$ is a $\g$-homotopy operator with respect to the
contraction operators $\iota_\xi^0$. Hence
$$ \iota_\xi^1(\psi)=(\iota_\xi^0((c \cdot \phi^{-1})\circ h)) \cdot \phi
-((c \cdot \phi^{-1})\circ h)\cdot \iota_\xi^1\phi=0.$$
On the other hand, using $c\circ i_{W\g}=0$ we have
$$ \d(((c \cdot \phi^{-1})\circ h)=(c \cdot \phi^{-1})\circ \d(h)
=(c \cdot \phi^{-1})\circ(\on{id}-i\circ \pi)
=c \cdot \phi^{-1}
$$
and therefore $\d(\psi)=c$.  Since obviously $L_\xi(\psi)=0$, the
proof is complete.
\end{proof}

\begin{remark}  \label{rem:psi1}
The homotopy $\psi$ constructed above vanishes on the unit of
$W\g$. This property is implied by $h(1)=0$ and $\Delta(1)=1\otimes
1$. \end{remark}

\subsection{$\g$-differential algebras of Weil type}
In this Section we will prove a more precise version of Corollary
\ref{cor:algebra}. In particular, we will describe a class of locally
free $\gda$'s of {\em Weil type} for which the algebra homomorphism
\eqref{eq:alghom} is in fact an isomorphism.  The results will not be
needed for most of our applications, except in Sections
\ref{subsec:isotropic} and \ref{sec:universal}.
\begin{definition}
Let $W$ be a locally free $\gda$, together with a $\ds$ homomorphism
$\pi:\,W\to \F$ such that $\pi\circ i=\on{id}$, where $i:\,\F\to W$ is
the unit for $W$.  Then $W$ will be called {\em of Weil type} if there
exists a homotopy operator $h$ between $i\circ \pi$ and $\on{id}$,
with $h\circ i=0$, such that $[L_\xi,h]=0$ and such that 
$h$ and all $\iota_\xi$ have degree $<0$ with respect to
some filtration
$$W=\bigcup_{N\ge 0} W^{(N)},\ W^{(0)}\subset W^{(1)}\subset \cdots$$
of $W$. 
\end{definition}
The definition is motivated by ideas from Guillemin-Sternberg
\cite[Section 4.3]{gu:su}.  In most examples, the $\Z_2$-grading on
$W$ is induced from a $\Z$-grading $W=\bigoplus_{n\ge 0} W^n$,
$\F\subset W^0$, with $\d$ of degree $+1$, $L_\xi$ of degree $0$ and
and $\iota_\xi$ of degree $-1$. If a homotopy operator $h$ with
$[L_\xi,h]=0$ exists, its part of degree $-1$ is still a homotopy
operator, and has the required properties.  (Note that compatibility
of the grading with the algebra structure is not needed.)  In
particular, the Weil algebra $W\g$ and its non-commutative version
$\ti{W}\g$ from Remark \ref{rem:wg1} are of Weil
type. The tensor product of two $\gda$'s of Weil type is again of 
Weil type. 
In the following Sections we will encounter several other
examples.

\begin{theorem}\label{th:technical}
Suppose $W,W'$ are $\gda$'s of Weil type. Then there exists a
$\gds$ homomorphism $\phi:\,W\to W'$ taking the unit in $W$ to the 
unit in $W'$. Moreover, any such $\phi$ is a $\g$-homotopy equivalence. 
\end{theorem}

It follows that all $\gda$'s of Weil type are $\g$-homotopy 
equivalent. The proof of Theorem \ref{th:technical} is somewhat
technical, and is therefore deferred to the appendix.

By Theorem \ref{th:technical}, many of the usual properties 
of $W\g$ extend to $\gda$'s of Weil type. For instance, since the 
basic cohomology of $W\g$ is $(S\g^*)_{\on{inv}}$, the same is 
true for any $\gda$ of Weil type: 

\begin{corollary}
If $\A=W$ is a $\gda$ of Weil type the Chern-Weil map 
$(S\g^*)_{\on{inv}}\to H(\A_{\on{basic}})$ 
is  an algebra isomorphism. 
\end{corollary}

\begin{corollary}
Let $W$ be a $\gda$ of Weil type and $\A$ be a locally free
$\gda$. 
\begin{enumerate}
\item
There is a $\gds$ homomorphism $\phi: W\to \A$ 
taking the unit in $W$ to the unit in $\A$. The induced map 
in basic cohomology is an algebra homomorphism.
\item
Any two $\gds$ homomorphisms $\phi: W\to \A$ that agree on the 
unit element of $W$ are $\g$-homotopic.  
\end{enumerate}
\end{corollary}
\begin{proof}
Both facts have already been established for $W=W\g$. Hence, by 
Theorem \ref{th:technical} they extend to arbitrary $\gda$'s of Weil type. 
\end{proof}
\begin{corollary}
The symmetrization map $W\g\to \ti{W}\g$ is a $\g$-homotopy equivalence, 
with $\g$-homotopy inverse the quotient map 
$\ti{W}\g\to W\g$. 
\end{corollary}
\begin{proof}
The composition $\ti{W}\g\to W\g\to \ti{W}\g$ is a $\gds$ homomorphism 
taking units to units, hence is $\g$-homotopic to the identity map. 
\end{proof}
\begin{corollary}\label{cor:B}
Let $W$ be of Weil type, and $\A$ a locally free $\gda$. 
Let $\phi:\,W\to \A$ be a $\gds$ homomorphism taking units to units. 
Then $\phi$ induces an algebra homomorphism in basic cohomology. 
More precisely, it intertwines the products
$m_W:\,W\otimes W\to W$ and 
$m_\A:\,\A\otimes\A\to \A$ 
up to $\g$-chain homotopy. 
\end{corollary}
\begin{proof}
The assertion is that the two maps $W\otimes W\to \A$ given by $m_\A
\circ (\phi \otimes \phi)$ and $\phi \circ m_{W}$ are $\g$-homotopic. 
But this follows
because $W\otimes W$ is a $\gda$ of Weil type, and each of the two
maps are $\gds$ homomorphisms that agree on the unit.
\end{proof}

\begin{remark}
The statement of Corollary \ref{cor:B} may be strengthened to the fact
that $\phi$ is an $A_\infty$-{\em morphism} ( see e.g. \cite{ma:ho}
for a precise definition).  Thus $\phi=\phi^{(1)}$ is the first term
in a sequence of $\wt{\g}$-equivariant maps $\phi^{(n)}:\,W^{\otimes
n}\to \A$. The second term $\phi^{(2)}=\psi$ is the $\g$-homotopy
between the two maps $W\otimes W\to \A$ given by $m_\A\otimes
(\phi\otimes\phi)$ and $\phi\otimes m_W$.  To construct $\phi^{(3)}$,
note that the two maps $W^{\otimes 3}\to \A$,
$$ m_\A \circ(\psi \otimes \phi) + \psi \circ (m_{W} \otimes
1),\ \ \ \ m_\A(\phi \otimes \psi) + \psi(1 \otimes m_W) $$
are $\gds$-homomorphisms. Since $\psi$ vanishes on the unit of 
$W^{\otimes 2}$ (cf. Remark \ref{rem:psi1}), both of these maps 
vanish at the unit of $W^{\otimes 3}$. 
It follows that there exists a $\g$-homotopy
$\phi^{(3)}:\,W^{\otimes 3}\to \A$ between these two maps, which 
again vanishes at the unit. Proceeding in this manner, one constructs 
the higher $\g$-homotopies. 
\end{remark}

\section{The Weil algebra $\ca{W}\g$.}\label{sec:caweil}
In this Section we construct, for any {\em quadratic} Lie algebra 
$\g$, an interesting non-commutative $\gda$ of Weil type.  

\subsection{Quadratic Lie algebras}
We begin by recalling some examples and facts about quadratic Lie 
algebras. From now on, we will refer to any non-degenerate symmetric bilinear
form on a vector space as a {\em scalar product}. A Lie algebra $\g$
with invariant scalar product $B$ will be called a {\em quadratic Lie
algebra}. First examples of quadratic Lie algebras are semi-simple 
Lie algebras, with $B$ the Killing form. Here are some other examples: 

\begin{examples}\label{ex:quadr}
\begin{enumerate}
\item
Let $\g$ be any Lie algebra, with given symmetric bilinear form.
Then the radical $\mf{r}$ of the bilinear form is an ideal, and the 
quotient  $\g/\mf{r}$ with induced bilinear form is quadratic. 
\item\label{ex:heis}
Let $\F^{2n}$ be equipped with the standard symplectic form,
$\om(e_{2i-1},e_{2i})=1$. Recall that the Heisenberg Lie algebra $H_n$
is the central extension 
$$ 0\lra \F \lra H_n \lra \F^{2n}$$
of the Abelian Lie algebra $\F^{2n}$ by $\F$,
with bracket defined by the cocycle $\om$. Let 
$c$ denote the basis vector for the center $\F\subset H_n$. 
Let another copy of $\F$, with basis vector $r$, act
on $\F^{2n}$ by infinitesimal rotation in each
$e_{2i-1}$--$e_{2i}$-plane: $r.e_{2i-1}=e_{2i},\ r.e_{2i}=-e_{2i-1}$.
This action lifts to derivations of $H_n$, and we may form the 
semi-direct product
$$ \g=\F\ltimes H_n.$$
The Lie algebra $\g$ is quadratic, with bilinear form given 
by $B(e_{2i-1},e_{2i})=B(c,r)=1$, and all other scalar 
products between basis vectors equal to $0$.
\item \label{ex:semidirect}
Let $\mf{s}$ be any Lie algebra, acting on its dual by the 
co-adjoint action. Viewing $\mf{s}^*$ as an Abelian Lie algebra, 
form the semi-direct product $\g=\mf{s}\ltimes\mf{s}^*$. 
Then $\g$ is a quadratic Lie algebra, with bilinear form  
$B$ given by the natural pairing between $\mf{s}$ and $\mf{s^*}$. 
More generally, given an invariant element 
$C\in (\wedge^3\mf{s})_{\mf{s}-\on{inv}}$, one obtains a quadratic 
Lie algebra where the bracket between elements of $\mf{s}^*$ is given 
by $[\mu,\mu']_\g=C(\mu,\mu',\cdot)\in \mf{s}$. See Section 
\ref{sec:examples} below. 
\end{enumerate}
\end{examples}
{} From a given quadratic Lie algebra $(\mf{a},B_{\mf{a}})$, new examples
are obtained by the {\em double extension construction} of
Medina-Revoy \cite{med:ca}: Suppose a second Lie algebra 
$\mf{s}$ acts on  $\mf{a}$ by derivations preserving the scalar product. 
Let $\om$ be the following $\mf{s}^*$-valued cocycle on $\mf{a}$,
$$ \l\om(a_1,a_2),\xi\r=B_{\mf{a}}(a_1,\xi.a_2),\ \
a_i\in\mf{a},\,\xi\in\mf{s}$$
and $\mf{a}\oplus\mf{s}^*$ the central extension of $\mf{a}$ defined
by this cocycle. The Lie algebra $\mf{s}$ acts on
$\mf{a}\oplus\mf{s}^*$ by derivations, hence we may form the
semi-direct product $\g=\mf{s}\ltimes (\mf{a}\oplus\mf{s}^*)$. The
given scalar product on $\mf{a}$, together with the scalar product on
$\mf{s}\ltimes\mf{s}^*$ given by the pairing, define a scalar product
on $\g$, which is easily checked to be invariant.
Notice that Example \ref{ex:quadr}(\ref{ex:heis}) is a special case of
this construction.

\subsection{The algebra $\W\g$}\label{subsec:main}
Consider the following 
example of the double extension construction (extended to super 
Lie algebras in the obvious way). Suppose $\g$ is any 
Lie algebra with a (possibly degenerate) symmetric bilinear form $B$. 
Then the super Lie algebra $\wt{\g}$ inherits an odd (!) 
symmetric bilinear form,  
$$ B_{\wt{\g}}(\xi,\xi')=0,\ B_{\wt{\g}}(\ol{\xi},\xi')=B(\xi,\xi'),\ 
B_{\wt{\g}}(\ol{\xi},\ol{\xi}')=0,\ \ \ \xi,\xi'\in\g.$$
The action of $\mf{s}=\F\d$ given by the differential on $\wt{\g}$
preserves $B_{\wt{\g}}$. The corresponding cocycle
$\om:\,\wt{\g}\otimes\wt{\g}\to \F$ is given by $B$ on the odd part
$\wt{\g}^\1=\g$ and vanishes on the even part. Thus, we obtain a
central extension $\wt{\g}\oplus \F\mf{c}$ by an even generator
$\mf{c}$ dual to $\d$; the new brackets between odd generators read,
$$ [\zeta,\zeta']_{\wt{\g}\oplus \F\mf{c}}=B(\zeta,\zeta')\mf{c}$$
while the brackets between even generators or between even and odd
generators are unchanged. The second step of the double extension
constructs the super Lie algebra 
\begin{equation}\label{eq:double}
\F\d\ltimes(\wt{\g}\oplus \F\mf{c}),
\end{equation} 
together with an odd symmetric bilinear form. The latter is 
non-degenerate if and only if $B$ is non-degenerate.

The super Lie algebra \eqref{eq:double} is a $\gdl$, where 
$\F\d\ltimes \wt{\g}$ acts by inner derivations. It contains the central 
extension $\wt{\g}\oplus \F\mf{c}$ as a $\g$-differential Lie subalgebra. 
Explicitly, the $\gds$ structure on $\wt{\g}\oplus \F\mf{c}$
is given by 
$$
L_\xi\ol{\zeta}=\ol{[\xi,\zeta]_\g},\ \
L_\xi\zeta=[\xi,\zeta]_\g,\ \
\iota_\xi\ol{\zeta}=[\xi,\zeta]_\g,\ \ 
\iota_\xi\zeta=B(\xi,\zeta)\mf{c},\ \ 
$$
while $\iota_\xi,L_\xi$ vanish on $\mf{c}$. 

\begin{remark}
It is instructive to compare the definition of the super Lie algebra
\eqref{eq:double} to the standard construction of affine Lie
algebras. Let $\g$ be a Lie algebra with invariant symmetric bilinear
form $B$.  Tensoring with Laurent polynomials, define an infinite
dimensional Lie algebra $\g[z, z^{-1}]=\g \otimes \F[z, z^{-1}]$ with
bilinear form $B'(x_1\otimes f_1,x_2\otimes f_2)= B(x_1,x_2)
\on{Res}(f_1f_2)$ where the residue ${\rm Res}$ picks the coefficient
of $z^{-1}$. The derivation $\partial(x\otimes f)=x\otimes \partial
f/\partial z$ preserves the bilinear form since ${\rm Res}(\partial
f/\partial z)=0$.  The double extension of $\g[z, z^{-1}]$
with respect to the derivation $\partial$ is called an affine Lie
algebra (at least if $B$ is non-degenerate). In a similar fashion,
letting $u$ be an {\em odd} variable we may tensor with the 
algebra $\F[u]=\{a+bu|\,a,b\in\F\}$ to define $\g[u]=\g\otimes 
\F[u]$. It carries an odd symmetric bilinear form, defined 
similar to $B'$ but with $\on{Res}$ replaced by the Berezin integral 
$\on{Ber}(a+bu)=b$. Again, the derivation $\d(x\otimes
f)=x\otimes \partial f/\partial u$ preserves the inner product since
$\on{Ber}(\partial f/\partial u)=0$.  Then $\g[u]\cong \wt{\g}$, the
derivation $\d$ is the Koszul differential, and the double extension
yields \eqref{eq:double}.

The two constructions may be unified to define a super Lie algebra
$\g[z,z^{-1},u]$ known as {\em super-affinization}, see \cite{ka:ve}. 
\end{remark}

We define the $\gda$ $\W\g$ as a quotient of the enveloping algebra, 
\begin{equation}\label{eq:weyl}
 \ca{W}\g:= U(\wt{\g}\oplus\F\mf{c})/<\mf{c}-1>.
\end{equation}
Note that $\W\g$ can be defined directly as a quotient of the
tensor algebra $\ca{T}(\wt{\g})$, in terms of generators $\xi,\ol{\xi}$
($\xi\in\g$) with relations $[\ol{\xi},\ol{\xi}']=\ol{[\xi,\xi']_\g}$,
$[\ol{\xi},{\xi}']=[\xi,\xi']_\g$ and $[\xi,\xi']=B(\xi,\xi')$.  In
particular there is a symmetrization map $S(\wt{\g})\to \W\g$. Let
$S(\wt{\g})$ carry the structure of a $\gda$, induced by its
identification with $S(\wt{\g}\oplus\F\mf{c})/<\mf{c}-1>$.
\begin{lemma}
The symmetrization map 
\begin{equation}\label{eq:symm}
 \ca{Q}_\g:\,\,S(\wt{\g})\to \W\g
\end{equation}
is a $\gds$ isomorphism. 
\end{lemma}
\begin{proof}
By Lemma \ref{lem:gds}, the PBW isomorphism 
$S(\wt{\g}\oplus \F\mf{c})\to U(\wt{\g}\oplus \F\mf{c})$ is a 
$\gds$ isomorphism. After quotienting the ideal $<\mf{c}-1>$ 
on both sides, the Lemma follows. 
\end{proof}

\begin{proposition}\label{prop:B=0}
If $B=0$, the inclusion map $\F\hra \W\g$ is a $\g$-homotopy equivalence.
\end{proposition}
\begin{proof}
It suffices to show that the inclusion $\F\hra S(\wt{\g})$ is a
$\g$-homotopy equivalence. Since $B=0$, $\wt{\g}$ (not just $\wt{\g}\oplus
\F\mf{c}$) is a $\gds$, and the map $s:\,\wt{\g}\to \wt{\g}$ (cf. the end of
Section \ref{subsec:ho}) entering the definition of Koszul homotopy
operator $h:\,\,S(\wt{\g})\to S(\wt{\g})$ commutes with Lie derivatives and
contractions. It follows that $h$ is a $\g$-homotopy in this case.
\end{proof}

\subsection{The non-degenerate case}
Let us now assume that $\g$ is a quadratic Lie algebra, i.e. that the
bilinear $B$ on $\g$ is {\em non-degenerate}. In this case, $\W\g$ 
becomes a $\gda$ of Weil type.  Indeed, the scalar
product $B$ defines an isomorphism $B^\sharp:\,\g^*\to\g$, and hence
$E_{\g^*}\cong E_\g=\wt{\g}$. This isomorphism identifies the $\gds$
structures on $E_{\g^*}\oplus \F\mf{c}$ and $\wt{\g}\oplus \F\mf{c}$,
hence it identifies $S(\wt{\g})$ with the Weil algebra
$W\g=S(E_{\g^*})$.  That is, for any invariant scalar product $B$,
\eqref{eq:symm} becomes a $\gds$ isomorphism
\begin{equation}\label{eq:Q}
 \ca{Q}_\g:\ W\g\to \ca{W}\g.
\end{equation}
Since $W\g$ is of Weil type, so is $\ca{W}\g$.  We will refer to
$\ca{Q}_\g$ as the {\em quantization map}.
\begin{remark}
On the subalgebras $\wedge\g\subset W\g$ resp.  $\on{Cl}(\g)\subset
\W\g$ generated by odd elements $\zeta$, the quantization map 
restricts to the usual quantization map (Chevalley symmetrization map) 
$$ q:\,\wedge\g\to \on{Cl}(\g)$$ 
for Clifford algebras, while on the subalgebras generated by even elements 
$\ol{\zeta}$ it becomes the PBW symmetrization map $S\g\to U\g$. 
\end{remark}
\begin{remark}
The Weil algebra $\ca{W}\g$ carries a connection $\g^*\cong\g \to
\ca{W}\g$, induced from the inclusion of 
$\g=\wt{\g}^\1$.  By construction, the quantization map \eqref{eq:Q} is 
the characteristic homomorphism for this connection. 
\end{remark}

As mentioned above, the horizontal subalgebra $(W\g)_{\on{hor}}$ of
the Weil algebra $W\g$ is isomorphic to the symmetric algebra, and the
differential $\d$ vanishes on $(W\g)_{\on{basic}}=(S\g^*)_{\on{inv}}$.
We will now show that similarly, the horizontal subalgebra
$(\ca{W}\g)_{\on{hor}}$ is $\g$-equivariantly isomorphic to the
enveloping algebra $U\g$. Let $\gamma:\,\g\to \on{Cl}(\g)$ be the map,
$$\gamma(\zeta)=\hh \sum_a [\zeta,e_a]_\g e^a$$
where $e_a$ is a basis of $\g$ and $e^a$ the dual basis with respect
to the given scalar product. It is a standard fact that the map 
$\gamma$ is a Lie algebra homomorphism, and  that 
Clifford commutator with $\gamma(\zeta)$ 
is the generator for the adjoint action of $\zeta$ 
on the Clifford algebra: 
$$ [\gamma(\zeta),\cdot]=L_\zeta^{\on{Cl}(\g)}$$
where the bracket denotes the super commutator in the Clifford algebra. 
Recall the definition of the map
$\lambda:\,\g^*\to \wedge^2\g^*$ by
$\iota_\xi\lambda(\mu)=-\ad_\xi^*\mu$. Identify $\g^*\cong\g$
by means of the scalar product, and let $q:\,\wedge\g\to \on{Cl}(\g)$ the
quantization map (i.e. symmetrization map) for the Clifford algebra. 
Then 
$$ \gamma(\zeta)=q(\lambda(\zeta)).$$
The curvature of the canonical connection on $\W\g$ is the 
map 
$$ \g\to (\W\g)_{\on{hor}},\ \ \zeta\mapsto \wh{\zeta}=\ol{\zeta}-\gamma(\zeta).$$
\begin{theorem}\label{prop:factor}\label{th:algebraiso}
The super algebra $\W\g$ is a tensor product
\begin{equation}\label{eq:tensorproduct} 
\W\g=U\g\otimes\on{Cl}(\g)
\end{equation}
where $U\g$ is generated by the even variables $\wh{\zeta}$ and the
Clifford algebra $\on{Cl}(\g)$ is generated by the odd variables
$\zeta$. Under this identification, the map 
$\ca{Q}_\g:\,W\g\to \W\g$ restricts to a vector space isomorphism
\begin{equation}\label{eq:dufsym}
S\g=(W\g)_{\on{hor}}\to U\g=(\W\g)_{\on{hor}}.
\end{equation}
In fact, \eqref{eq:dufsym} is an algebra isomorphism on $\g$-invariants.
\end{theorem}
\begin{proof}
The elements $\wh{\zeta}\in\W\g$ are the images of the corresponding elements 
in $W\g$ (denoted by the same symbol) under the quantization map 
$\ca{Q}_\g$. 
The commutator of two such elements $\wh{\zeta},\wh{\zeta'}\in\W\g$ 
is given by 
$$ [\wh{\zeta},\wh{\zeta}']=\ol{[\zeta,\zeta']_\g}
-L_\zeta\gamma(\zeta')-L_{\zeta'}\gamma(\zeta)+\gamma([\zeta,\zeta']_\g)
=\wh{[\zeta,\zeta']_\g}.$$
Hence the variables $\wh{\zeta}$ generate a copy of the enveloping
algebra $U\g\subset \W\g$. On the other hand, the odd variables
$\zeta$ generate a copy of the Clifford algebra.  Since
$[\zeta,\wh{\zeta'}]=0$ for $\zeta,\zeta'\in\g$, the decomposition
\eqref{eq:tensorproduct} follows. 
Since the map $\ca{Q}_\g:\,W\g\to \ca{W}\g$ is a $\gds$ isomorphism,
and the Weil differential $\d$ vanishes on basic elements 
in $W\g$, it also vanishes on $(\ca{W}\g)_{\on{basic}}$.
Hence $H((\ca{W}\g)_{\on{basic}})=(\ca{W}\g)_{\on{basic}}=
(U\g)_{\on{inv}}$, and the last claim of the Theorem follows from Corollary 
\ref{cor:algebra}.
\end{proof}

\section{Duflo isomorphism}\label{sec:duflo}
We will now show that the isomorphism
\eqref{eq:dufsym} from $S\g$ to $U\g$
is exactly the Duflo isomorphism, for the case of
quadratic Lie algebras. Thus Theorem \ref{th:algebraiso} proves
Duflo's theorem for this case.  More generally, we will show that the
map $\ca{Q}_\g:\,W\g\to \W\g$ coincides with the quantization map
introduced in \cite{al:no}.
\begin{proposition}
Each of the derivations $\iota_\xi,\d,L_\xi$ of $\W\g$ is inner. 
\end{proposition} 
\begin{proof}
By construction, $\iota_\xi=[\xi,\cdot]$ and
$L_\xi=[\ol{\xi},\cdot].$ To show that the differential $\d$ is inner,
choose a basis $e_a$ of $\g$, and let $e^a$ be the dual basis with
respect to $B$. Then ${\sum}_a\ol{e_a}e^a\in\W\g$ is an invariant
element, independent of the choice of basis.  We have
$$ [{\sum}_a\ol{e_a}e^a,\ol{\xi}]=0,\ \
[{\sum}_a\ol{e_a}e^a,\xi]=\ol{\xi}-2\gamma(\xi).$$
Let 
$$ \ca{D}:={\sum}_a\ol{e_a}e^a-\f{2}{3}{\sum}_a \gamma(e_a)e^a.$$
Since $\ca{D}$ is invariant, $[\ca{D},\ol{\xi}]=-L_\xi \ca{D}=0=\d\ol{\xi}$. 
The element $\phi=\f{1}{3}{\sum}_a \gamma(e_a)e^a$ satisfies 
$[\phi,\xi]=\gamma(\xi)$. Therefore
On the other hand, since $[\f{1}{3}{\sum}_a \gamma(e_a)e^a,\xi]=\gamma(\xi)$, 
we have 
$$ [\ca{D},\xi]=\ol{\xi}+\sum_a \ad_\xi(e_a) e^a -2\gamma(\xi)=\ol{\xi}
=\d\xi.$$ 
\end{proof}

\begin{remarks}\label{rem:D}
\begin{enumerate}
\item\label{rem:D1}
The cubic element $\ca{D}$ may be interpreted as a quantized chain of 
transgression. Indeed, it is easily checked that $\ca{Q}_\g(D)=\ca{D}$
where 
$$ D:=h(\sum_a \wh{e_a} \wh{e^a})={\sum}_a\ol{e_a}e^a-\f{2}{3}\sum_a \lambda(e_a) e^a\in W\g$$
is the chain of transgression corresponding to the 
quadratic polynomial 
$\sum_a \wh{e_a} \wh{e^a}\in (S\g)_{\on{inv}}\subset W\g$. 
Here  $h$ is the standard homotopy operator for the Weil algebra. 
\item
The fact that the derivation $\d$ is inner may also be formulated in
terms of the quadratic super Lie algebra $\F\d\ltimes(\wt{\g}\oplus
\F\mf{c})$. Indeed, it may be verified that the cubic element
$$\mf{c}^2 \d - \mf{c} {\sum}_a\ol{e_a}e^a + \f{2}{3}{\sum}_a \gamma(e_a)e^a
$$
in the enveloping algebra of $\F\d\ltimes(\wt{\g}\oplus \F\mf{c})$ is a
central element. Specializing to $\mf{c}=1$ we see that $\d-D$ is
central in the quotient by $< \mf{c}- 1 >$.
\end{enumerate}
\end{remarks}

We now recall the definition of the Duflo map $S\g\to U\g$. Let $
\ol{S\g^*}=\coprod_{k=0}^\infty S^k\g^*$ be the completion of the
symmetric algebra, or equivalently  the algebraic dual space
to $S\g$. Informally, we will view $\ol{S\g^*}$ as Taylor series 
expansions of functions on $\g$. There is an algebra
homomorphism
$$\ol{S\g^*}\to \End(S\g),\ F\mapsto \wh{F}$$ 
extending the natural action of $\g^*$ by derivations, that is, 
$\hat{F}$ is an 
infinite order differential operator acting on polynomials. 
Let $j(z)=\f{\sinh(z/2)}{z/2}$ and define 
$J\in \ol{S\g^*}$ by
$$ J(\xi)=\det(j(\ad_\xi))=e^{\on{tr}\big((\ln j)(\ad_\xi)\big)}.$$
The square root of $J$ is a well-defined element of $\ol{S\g^*}$. 
The Duflo map is the composition
\begin{equation}\label{eq:duflo}
 \on{sym}_{U\g}\circ \wh{J^{1/2}}:\,S\g\to U\g.
\end{equation}
The quantization map in \cite{al:no} is an extension of the Duflo
map, for the case that $\g$ is quadratic. Write $W\g=S\g\otimes\wedge\g$, as in
\eqref{eq:factWg}, and $\W\g=U\g\otimes \on{Cl}(\g)$ as in
\eqref{eq:tensorproduct}.  
Let 
$$(\ln j)'(z)=\f{1}{2}\coth\f{z}{2}-\f{1}{z}$$
be the logarithmic
derivative of the function $j$, and let 
$\mf{r} \in \ol{S\g^*} \otimes \wedge^2\g$ be given by 
$r(\xi)=(\on{ln}j)'(\ad_\xi)$, where we identify skew-symmetric operators on
$\g$ with elements in $\wedge^2\g$. Put 
$$ \ca{\S}(\xi)=J^{1/2}(\xi)\exp(\mf{r}(\xi))$$ 
and let $ \wh{\iota(\ca{S})}\in
\on{End}(W\g)$ denote the corresponding operator, where the
$\ol{S\g^*}$ factor acts as an infinite order differential operator
on $S\g$, and the $\wedge\g$ factor acts by contraction on $\wedge\g$.
Let $q:\,\wedge\g\to \on{Cl}(\g)$ be the Chevalley quantization map
for the Clifford algebra. The tensor product of the PBW symmetrization
map $\on{sym}_{U\g}:\,S\g\to U\g$ and the Chevalley quantization map
$q:\,\wedge\g\to \on{Cl}(\g)$ define a linear isomorphism
$\on{sym}_{U\g}\otimes q:\,W\g\to \W\g$.  Put differently, this is the
symmetrization map with respect to the generators $\xi,\wh{\xi}$,
rather than the generators $\xi,\ol{\xi}$ used in the definition
of $\ca{Q}_\g$.
\begin{theorem}\label{th:same}
Under the identification $W\g=S\g\otimes\wedge\g$ and 
$\W\g=U\g\otimes\Cl(\g)$, the quantization map 
is given by the formula, 
$$\ca{Q}_\g=(  \on{sym}_{U\g}\otimes\, q)\circ \wh{\iota(\ca{S})}:\ 
W\g\to\W\g.
$$
In particular, its restriction to the symmetric algebra $S\g$ is the
Duflo map.
\end{theorem}
\begin{proof}
We use an alternative description of the symmetrization map $S(\wt{\g})\to
\ca{T}(\wt{\g})$.  Let $\nu^a\in E_{\g^*}^\1$
and $\mu^a\in E_{\g^*}^\0$ be ``parameters''. The symmetrization
map is characterized by its property that for all $p$, the map
$\on{id}\otimes \ca{Q}_\g$ takes 
the $p$th power of 
${\sum}_a(\nu^a e_a+\mu^a \ol{e_a})$ 
in the algebra $S(E_{\g^*})\otimes W\g$
to the corresponding $p$th power in 
$S(E_{\g^*})\otimes\W\g$. These conditions may be combined into a 
single condition  
$$\on{id}\otimes \ca{Q}_\g:\, \exp_{W\g}({\sum}_a(\nu^a e_a+\mu^a
\ol{e_a})) \mapsto \exp_{\W\g}({\sum}_a(\nu^a e_a+\mu^a \ol{e_a}));$$
here the exponentials are well-defined in completions 
$\ol{S(E_{\g^*})}\otimes W\g$ and $\ol{S(E_{\g^*})}\otimes\W\g$, respectively.

We want to re-express the symmetrization map in terms of the 
new generators $e_a,\wh{e_a}=\ol{e_a}-\lambda(e_a)$ of $W\g$ and 
$e_a,\wh{e_a}=\ol{e_a}-\gamma(e_a)$ of $\W\g$. Using that $e_a$ and $\wh{e_a}$ 
commute in $\W\g$, we may separate the $\Cl(\g)$ and $U\g$-variables in 
the exponential and obtain: 
$$ \exp_{\W\g}({\sum}_a(\nu^a e_a+\mu^a \ol{e_a}))=
\exp_{\Cl(\g)}\big({\sum}_a(\nu^a e_a+\mu^a \gamma(e_a)\big)
\exp_{U\g}({\sum}_a\mu^a \wh{e_a}).
$$ 
The factor $\exp_{U\g}({\sum}_a\mu^a \wh{e_a})$ is the image of
$\exp_{S\g}({\sum}_a\mu^a \wh{e_a})$ under the symmetrization map
$\on{sym}_{U\g}:\,S\g\to U\g$. The other factor is the exponential of
a quadratic expression in the Clifford algebra. Using \cite[Theorem
2.1]{al:cli} such exponentials may be expressed in terms of 
the corresponding exponentials in the exterior algebra:
$$ \exp_{\Cl(\g)}{\sum}_a\big(\nu^a e_a+\mu^a \gamma(e_a)\big)
=q\Big(\iota(\S(\mu))\exp_{\wedge\g}{\sum}_a\big(\nu^a e_a+\mu^a
\lambda(e_a)\big)\Big),
$$
where $\iota:\,\wedge\g\to \on{End}(\wedge\g)$ is contraction. 
This shows 
\beq 
\lefteqn{\exp_{\W\g}{\sum}_a(\nu^a e_a+\mu^a \ol{e_a})}\\
&=&(\on{sym}_{U\g}\otimes q)
\Big(\iota(\S(\mu))\exp_{W\g}{\sum}_a\big(\nu^a e_a+\mu^a
(\wh{e_a}+\lambda(e_a)) \big)\Big)\\
&=& \ca{Q}_\g\Big(\exp_{W\g}{\sum}_a\big(\nu^a e_a+\mu^a(\wh{e_a}+\lambda(e_a)\big)\Big)\\
&=& \ca{Q}_\g\Big(\exp_{W\g}{\sum}_a\big(\nu^a e_a+\mu^a \ol{e_a}\big)\Big). 
\eeq
\end{proof}
A different proof of Theorem \ref{th:same} will be given in Section
\ref{sec:rouviere}.  Our result
shows that while the Duflo map itself is not a symmetrization map, it
may be viewed as the restriction of a symmetrization map for a larger
algebra.
Using Theorem \ref{th:same}, we obtain a very simple proof of the following 
result from \cite{al:no,ko:cu}: 
\begin{proposition}\label{prop:square}
The square of the 
cubic element $\ca{D}$ is given by 
\begin{equation}\label{eq:casimir}
 \ca{D}^2=\f{1}{2} \on{Cas}_\g
+\f{1}{48} \on{tr}_\g(\on{Cas}_\g)
\end{equation}
where $\on{Cas}_\g=\sum_a \wh{e_a} \wh{e^a}\in (U\g)_{\on{inv}}\subset
\W\g$ is the quadratic Casimir element, and $\on{tr}_\g(\on{Cas}_\g)$
its trace in the adjoint representation.
\end{proposition}
\begin{proof}[Outline of proof]
Write $\ca{D}=\ca{Q}_\g(D)$ as in Remark \ref{rem:D} (\ref{rem:D1}).
Then 
$$ \ca{D}^2=\hh [\ca{D},\ca{D}]=
\hh \d(\ca{D})=\hh \ca{Q}_\g(\d D)=
\hh \ca{Q}_\g(\sum_a \wh{e_a} \wh{e^a}). $$
Using the explicit formula \eqref{eq:duflo} for the Duflo map, one
finds that $\ca{Q}_\g(\sum_a \wh{e_a}\wh{e^a})$ is equal to
$\on{Cas}_\g+\f{1}{24}\on{tr}_\g(\on{Cas}_\g)$.
\end{proof}

\begin{remarks}
\begin{enumerate}
\item The algebra $\W\g$ carries a natural $\Z$-filtration, where the
odd generators have degree $1$ and the even generators have degree
$2$. Its associated graded algebra is the Weil algebra $W\g$, with its
standard grading. The filtration on $\W\g$ is compatible with the
$\Z_2$-grading in the sense of \cite{ko:sy2}, and therefore induces
the structure of a graded Poisson algebra on $W\g$. On generators, the
formulas for the graded Poisson bracket are given by $\{\xi,\zeta\}=
B(\xi,\zeta)$, $\{\ol{\xi},\zeta\}=\ad_\xi\zeta$,
$\{\ol{\xi},\ol{\zeta}\}= \ol{\ad_\xi\zeta}$. The Weil differential
may be written as $\d=\{D,\cdot\}$ with $D$ as in Remark
\ref{rem:D} (\ref{rem:D1}).
\item
It is possible to re-introduce a grading on $\W\g$, by adding 
an extra parameter.
Let $\hbar$ be a variable of degree $2$, and view
$S(\F\cdot\hbar)$ as a $\gda$ with contractions, Lie
derivatives and differential all equal to zero. Define a $\gda$
$\ca{W}\g[\hbar]$ as the quotient of $\ca{T}(\wt{\g})\otimes
S(\F\cdot\hbar)$ by the relations,
$$ [\zeta,\zeta']=\hbar B(\zeta,\zeta'),\ \ 
[\ol{\zeta},\zeta']=\hbar [\zeta,\zeta']_\g,\ \ 
[\ol{\zeta},\ol{\zeta'}]=\hbar \ol{[\zeta,\zeta']_\g}
$$
(on the left hand side, the bracket denotes super commutators).  
Note that the ideal generated by these relations is graded. Hence
$\ca{W}\g[\hbar]$ is graded and the symmetrization map
$W\g[\hbar]=W\g\otimes S(\F\cdot\hbar)\to \ca{W}\g[\hbar]$ preserves
degrees.  The Weil algebra $W\g$ is obtained by dividing out the ideal
$\l\hbar\r$, while $\W\g$ is obtained by dividing out
$\l\hbar-1\r$. This clearly exhibits $\W\g$ as a deformation of the
Weil algebra $W\g$.
\item
Pavol Severa explained to us in the summer of 2001, that the quantization
map from \cite{al:no} is closely related to the exponential
map for a central extension of the super Lie group $T[1]G$, 
as discussed in his paper \cite{sev:so}. Theorem \ref{th:same} may be 
viewed as the algebraic version of Severa's observation.
\end{enumerate}  
\end{remarks}

\section{Vogan conjecture}\label{sec:qu}
Suppose $\k\subset \g$ is a Lie subalgebra admitting a $\k$-invariant 
complement $\mf{p}$. Thus $\g=\k\oplus \mf{p}$ with
$$ \ \ [\k,\k]_\g \subset \k,\ \ [\k,\mf{p}]_\g \subset \mf{p}.$$
Any $\gds$ $E$ becomes a $\kds$ by restricting the action to
$\wt{\k}\subset \wt{\g}$. If $\A$ is a $\gda$ with connection
$\theta:\,\g^*\to \A^\1$, then the restriction of $\theta$ to
$\k^*\subset \g^*$ defines a connection for $\A$, viewed as a
$\kda$. In particular, any $\gda$ of Weil type becomes a $\kda$ of
Weil type. Theorem \ref{th:technical} shows that the 
projection map $W\g\to W\k$ induced by the projection 
is a $\k$-homotopy equivalence, with homotopy inverse
induced by the inclusion $E_{\k^*}\to E_{\g^*}$. 
Suppose now that $\g$ is a quadratic Lie algebra, and 
that the restriction of the scalar product $B$ to the subalgebra 
$\k$ is again non-degenerate. We will refer to $\k$ as a {\em quadratic 
subalgebra}. In this case, we may take $\mf{p}$ to be the orthogonal 
complement of $\k$ in $\g$. 
\begin{example}
Let $\g$ be a semi-simple Lie algebra, with $B$ the Killing form,
and $\g=\k\oplus\mf{p}$ a Cartan decomposition. 
Then $\k$ and $\mf{p}$ are orthogonal, and $B$ is 
negative definite on $\k$ and positive definite on $\mf{p}$. 
See \cite[Chapter 3.7]{he:di} or \cite[Chapter VII.2]{kn:li}.  
\end{example}
\begin{example}
Suppose $\F=\C$. For any $\xi\in\g$, the generalized eigenspace 
for the $0$ eigenvalue of $\ad_\xi$
$$\k=\{\zeta\in\g|\,\ad_\xi^N\zeta=0\mbox{ for }N>>0\}$$ 
is a quadratic subalgebra, with $\mf{p}$ the direct sum of generalized
eigenspaces for nonzero eigenvalues. Indeed, 
given $\zeta_1\in\k,\ \zeta_2\in\mf{p}$, let $N>0$ with 
$\ad_\xi^N\zeta_1=0$. Since $\ad_\xi$ is invertible 
on $\mf{p}$, we have 
$$ B(\zeta_1,\zeta_2)=(-1)^N B(\ad_\xi^N \zeta_1,\ \ad_\xi^{-N}\zeta_2)=0.$$
This shows $\mf{p}= \k^\perp$. 
\end{example}
The inclusion
$\wt{\k}\oplus \F\mf{c}\to \wt{\g}\oplus \F\mf{c}$ is a $\kdl$ homomorphism, 
hence it extends to a $\kda$ homomorphism
$$ U(\wt{\k}\oplus \F\mf{c})\to U(\wt{\g}\oplus \F\mf{c}).$$
Taking quotients by the ideals generated by $\mf{c}-1$, we obtain a
$\kda$ homomorphism $\W\k\to \W\g$.  We
obtain a commutative diagram of $\kds$ homomorphisms,
\begin{equation}\label{eq:weilalg}
 \xymatrix{W\g\ar[r]_{\ca{Q}_\g} &\W\g\\ 
                      W\k\ar[r]_{\ca{Q}_\k}\ar[u]&\W\k\ar[u]},
\end{equation}
in which all maps are $\k$-homotopy equivalences, and the induced maps 
in basic cohomology are all algebra isomorphisms. We will now 
interpret these maps in terms of the isomorphism
$\W\g=U\g\otimes\Cl(\g)$.  We have
$$ (\W\g)_{\k-\on{hor}}=U\g\otimes \on{Cl}(\mf{p})$$
and therefore $ (\W\g)_{\k-\on{basic}}=(U\g\otimes 
\on{Cl}(\mf{p}))_{\k-\on{inv}}$. 
The Kostant cubic Dirac operator for the pair $\g,\k$ is 
defined to as the difference of the Dirac operators for 
$\g$ and $\k$: 
$$ \D_{\g,k}=\D_\g-\D_\k.$$
The first two parts of the following Proposition were proved by 
Kostant in \cite{ko:cu}. Let 
\begin{equation}
 \chi:\ U\k\hra U\g\otimes \on{Cl}(\mf{p})
\end{equation}
be the map given by the inclusion 
$(W\k)_{\k-\on{hor}}\hra (\W\g)_{\k-\on{hor}}$. 

\begin{proposition}
\begin{enumerate}
\item 
The cubic Dirac operator $\D_{\g,k}$ lies in the algebra $(U\g\otimes 
\on{Cl}(\mf{p}))_{\k-\on{inv}}$. 
\item 
The square of $\ca{D}_{\g,k}$ is given by the formula, 
$$ \D_{\g,\k}^2=\f{1}{2}\on{Cas}_\g-\f{1}{2}\chi(\on{Cas}_\k)
+\f{1}{48}(\on{tr}_\g(\on{Cas}_\g)-\on{tr}_\k(\on{Cas}_\k)).$$
\item
The restriction of the differential on $\W\g$ to the
subalgebra $(U\g\otimes \on{Cl}(\mf{p}))_{\k-\on{inv}}$ is 
a graded commutator $[\D_{\g,k},\cdot]$. 
\end{enumerate}
\end{proposition}
\begin{proof}
(a) It is clear that $\D_{\g,k}$ is $\k$-invariant. Furthermore, 
for $\xi\in \k$, we have 
$$[\D_{\g,k},\xi]=[\D_\g,\xi]-[\D_\k,\xi]=\ol{\xi}-\ol{\xi}=0$$
so $\D_{\g,\k}$ is $\k$-basic. 
(b) Since $\iota_\xi=[\xi,\cdot]$ and $L_\xi=[\ol{\xi},\cdot]$, the
basic subalgebra $(\W\g)_{\k-\on{basic}}$ is exactly the commutant of
the subalgebra $\W\k\subset\W\g$. Hence $\D_{\g,\k}\in
(\W\g)_{\k-\on{basic}}$ and $\D_\k\in\W\k$ commute, and the formula
follows from \eqref{eq:casimir} by squaring the identity
$\D_\g=\D_{\g,\k}+\D_\k$.
(c) On elements of $(\W\g)_{\k-\on{basic}}$, $[\D_\k,\cdot]$ vanishes
since $\D_\k\in \W\k$.  Hence $[\D_{\g,\k},\cdot]$ coincides with
$[\D_\g,\cdot]$ on $ (\W\g)_{\k-\on{basic}}$.
\end{proof}

The following theorem is a version of Vogan's conjecture (as formulated in 
Huang-Pandzic \cite{hu:vo}) for quadratic Lie 
algebras. It was first proved by Huang-Pandzic \cite[Theorems 3.4,
5.5]{hu:vo} for symmetric pairs, and by Kostant 
\cite[Theorem 0.2]{ko:di} for reductive pairs. In a recent paper, 
Kumar \cite{ku:ind} interpreted the Vogan conjecture in terms of 
induction maps in non-commutative equivariant cohomology.  

\begin{theorem}\label{th:vogan}
The map $\chi:\,(U\k)_{\k-\on{inv}}\to
(U\g\otimes\on{Cl}(\mf{p}))_{\k-\on{inv}}$ takes values in cocycles
for the differential $[\D_{\g,\k},\cdot]$, and descends to an algebra isomorphism from
$(U\k)_{\k-\on{inv}}$ to the cohomology of
$(U\g\otimes\on{Cl}(\mf{p}))_{\k-\on{inv}}$.
The map $(U\g)_{\g-\on{inv}}\to  (U\k)_{\k-\on{inv}}$
taking $z\in (U\g)_{\g-\on{inv}}$ to the cohomology class of 
$z\otimes 1\in (U\g\otimes\on{Cl}(\mf{p}))_{\k-\on{inv}}$ fits into 
a commutative diagram, 
$$ 
\xymatrix{(S\g)_{\g-\on{inv}}\ar[r]\ar[d]&(U\g)_{\g-\on{inv}}\ar[d]\\
(S\k)_{\k-\on{inv}}\ar[r]& (U\k)_{\k-\on{inv}}}
$$
where the horizontal maps are Duflo maps and the left vertical 
map is induced by the projection $\g\to \k$. 
\end{theorem}
\begin{proof}
The commutative diagram \eqref{eq:weilalg} gives rise to a 
commutative diagram of $\k$-basic subcomplexes, 
$$ 
\xymatrix{(S\k)_{\k-\on{inv}}\ar[r]\ar[d]&(U\k)_{\k-\on{inv}}\ar[d]_\chi\\
(S\g\otimes\wedge\mf{p})_{\k-\on{inv}}\ar[r]& 
(U\g\otimes\on{Cl}\mf{p})_{\k-\on{inv}}}
$$
As mentioned after \eqref{eq:weilalg}, all of the maps in this diagram 
induce algebra isomorphisms in cohomology. This proves the first 
part of the theorem. The second part follows by combining this diagram 
with a commutative diagram 
$$ 
\xymatrix{(S\g)_{\g-\on{inv}}\ar[r]\ar[d]&(U\g)_{\g-\on{inv}}\ar[d]\\
(S\g\otimes\wedge\mf{p})_{\k-\on{inv}}\ar[r]& 
(U\g\otimes\on{Cl}\mf{p})_{\k-\on{inv}}}
$$
given by the inclusion of $\g$-basic subcomplexes of $W\g,\W\g$ 
into the $\k$-basic subcomplexes.  
\end{proof}
\section{Harish-Chandra isomorphism}\label{sec:hc}
Let $\g$ be a quadratic Lie algebra, with scalar product $B$, and 
$\k\subset\g$ a quadratic subalgebra, with orthogonal complement 
$\mf{p}=\k^\perp$. In the previous Section, we obtained an
algebra homomorphism $(U\g)_{\g-\on{inv}}\to (U\k)_{\k-\on{inv}}$ 
which under the Duflo isomorphism corresponds to the natural 
projection $(S\g)_{\g-\on{inv}}\to (S\k)_{\k-\on{inv}}$. 
We will now describe an alternative construction of this map
for enveloping algebras, generalizing the Harish-Chandra construction 
\cite[Chapter V.5]{kn:li}

An extra ingredient for the Harish-Chandra 
map is a $\k$-invariant splitting $\mf{p}=\n_-\oplus \n_+$ into  
Lie subalgebras of $\g$ which are {\em isotropic}, i.e. 
such that $B$ vanishes on $\n_\pm$. Thus
\begin{equation}\label{eq:decomp} 
 \g=\n_-\oplus \k \oplus \n_+
\end{equation}
(direct sum of subspaces). 

\begin{examples}\label{ex:morex}
\begin{enumerate}
\item 
In the standard setting of the Harish-Chandra theorem, $\g$ is a
semi-simple Lie algebra over $\F=\C$, with compact real form $\g_\R$,
$\k=\t^\C$ is the complexification of a maximal Abelian subalgebra of
$\t\subset \g_\R$, and $\n_\pm$ are nilpotent subalgebras given as sums of root
spaces for the positive/negative roots. More generally, one could take
$\k$ to be the centralizer of some element $\xi\in\t$, and $\n_+$
(resp. $\n_-$) the direct sum of the positive (resp. negative) root
spaces that are not contained in $\k$.
\item\label{ex:euklid} Suppose $\F=\C$. Consider the
$2n+2$-dimensional nilpotent Lie algebra $\C\ltimes H_n$ from Example
\ref{ex:quadr}(\ref{ex:heis}).  We obtain a
decomposition \eqref{eq:decomp} by letting $\k$ be the Abelian
subalgebra spanned by $r$ and $c$, and letting $\n_\pm$ be the span of
$e_{2i}\pm \sqrt{-1} e_{2i-1}$.
\item \label{ex:eigenspace}
Both of these examples are special case of the following set-up.
Suppose that $\g_\R$ is an arbitrary quadratic Lie algebra over $\R$, 
$\g$ its complexification, and $\xi\in\g_\R$. Since 
$\ad_\xi$ preserves the quadratic form $B$, all
eigenvalues of $\f{1}{i}\ad_\xi$ are real. Let $\g_t\subset\g$ denote the
generalized eigenspace for the eigenvalue $t\in\R$. Then
$[\g_t,\g_{t'}]_\g \subset \g_{t+t'}$, and $B(\g_t,\g_{t'})=0$ for
$t+t'\not =0$. A decomposition \eqref{eq:decomp} is obtained by
setting $\k=\g_0,\ \n_-=\bigoplus_{t<0}\g_t,\
\n_+=\bigoplus_{t>0}\g_t$. 
\end{enumerate}
\end{examples}
By the Poincar\'{e}-Birkhoff-Witt theorem, the decomposition
\eqref{eq:decomp} of $\g$ yields a decomposition of the enveloping
algebra $U\g$,
$$ U\g=(\mf{n}_-\,U\g+U\g\,\mf{n}_+)\oplus U\k,$$
hence a (generalized) {\em Harish-Chandra projection}
$$\kappa_U:\,U\g\to U\k.$$ 
The projection $\kappa_U$ is $\k$-invariant, and restricts to an
algebra homomorphism on the subalgebra $\mf{n}_-\, U\g\, \mf{n}_+
\oplus U\k$.  Similar to $\kappa_U$, we define 
Harish-Chandra projections $\kappa_{\Cl}: \Cl(\g) \to
\Cl(\k)$ using the decomposition 
$$
\Cl(\g) = (\mf{n}_-\,\Cl(\g) +\Cl(\g) \,\mf{n}_+)\oplus \Cl(\k)
$$
and $\kappa_\W:\,\W\g\to \W\k$ using the decomposition 
\begin{equation}\label{eq:decomp2}
 \W\g=(\wt{\mf{n}_-}\W\g+ \W\g \,\wt{\mf{n}_+})\oplus \W\k.
\end{equation}
In Harish-Chandra's construction for enveloping algebras, it is necessary to 
compose the projection $\kappa_U$ with a ``shift''. Consider 
the infinitesimal character on $\k$, 
$$\zeta\mapsto \on{Tr}_{\n_+} \ad_\zeta.$$
The map 
$$\tau:\,\k\to U\k,\ \zeta\mapsto \zeta+\hh \on{Tr}_{\n_+} \ad_\zeta$$
is a Lie algebra homomorphism, hence it extends to an algebra automorphism, 
$\tau:\,U\k\to U\k$. In the standard case where $\k$ is a Cartan subalgebra
of a complex semi-simple Lie algebra, this is the ``$\rho$-shift''.
Remarkably, the shift is already built into the the projection 
$\kappa_\W$:  
\begin{proposition} \label{prop:kappas}
The following diagram commutes: 
$$
\xymatrix{\W\g\ar[r]_{\!\!\!\!\!\!\!\! \cong}
\ar[d]_{\kappa_{W}}&U\g\otimes\Cl(\g) \ar[d]^{(\tau \circ\kappa_U)\otimes \kappa_{\Cl} }\\ 
\W\k\ar[r]_{\!\!\!\!\!\!\!\! \cong} &U\k\otimes \Cl(\k)}
$$
\end{proposition}

\begin{proof}
Observe that the two projections $\kappa_\W$ and $\kappa_U \otimes
\kappa_\Cl$ both vanish on $E_{\mf{n}_-}\W\g+ \W\g
E_{\mf{n}_+}$. Hence it suffices to compare them on the subalgebra
$\W\k \subset \W\g$.  Elements in $\W\k$ may be written in the form
$\ol{\xi_{i_1}}\cdots \ol{\xi_{i_r}}\,x$ where
$\xi_{i_1},\ldots,\xi_{i_r}\in\k$ and $x\in \Cl(\k)\subset \W\k$.  For
$\xi\in\k$, we have $\ol{\xi}=\wh{\xi}+\gamma(\xi)$, where
$\gamma(\xi)$ decomposes into parts $\gamma^\k(\xi)\in \on{Cl}(\k)$
and $\gamma^{\mf{p}}(\xi)\in\Cl(\mf{p})$.  Let $b_i\in \mf{n}_-$ 
and $c_j\in \mf{n}_+$ be dual bases, i.e. $B(b_i,c_j)=\delta_{ij}$.
We have,
\beq
\gamma^{\mf{p}}(\xi)&=&\hh\sum_j\big((\ad_\xi b_j)c_j
+(\ad_\xi c_j)b_j\big)\\&=&
\hh\sum_j\big((\ad_\xi b_j)c_j-b_j(\ad_\xi c_j)
\big)-\hh \on{tr}_{\n_+}(\ad_{\xi}).
\eeq
Here we have used $ \sum_j [b_j,\ad_\xi c_j]=
\on{tr}_{\n_+}(\ad_\xi), $ for all $\xi\in\k$. 
Since $(\ad_\xi b_j)c_j-b_j(\ad_\xi c_j)\in \n_-\Cl(\g)\n_+$, 
it follows that 
$$ \ol{\xi}=\wh{\xi}-\hh \on{tr}_{\n_+}(\ad_{\xi})+\gamma^\k(\xi)\ \
\on{mod}\ \n_-\Cl(\g)\n_+.$$
Hence,
$$ \ol{\xi_{i_1}}\cdots \ol{\xi_{i_r}}\,x= \big((\wh{\xi}_{i_1}-\hh
\on{tr}_{\n_+}(\ad_{\xi_{i_1}})+\gamma^\k(\xi_{i_1}))\cdots
(\wh{\xi}_{i_r}-\hh \on{tr}_{\n_+}(\ad_{\xi_{i_r}})
+\gamma^\k(\xi_{i_r}))x\big)\ +\ldots$$
where the terms $\ldots$ are in $E_{n_-}\W\g +\W\g E_{n_+}$. 
The term in the large parentheses lies in the 
image of the tensor products of inclusions $U\k\hra U\g$, 
$\Cl(\k)\hra \Cl(\g)$, and $\kappa_U\otimes\kappa_{\Cl}$ 
is the identity map on this image. Comparing with 
$$\ol{\xi_{i_1}}\cdots \ol{\xi_{i_r}}\,x=
(\wh{\xi}_{i_1}+\gamma^\k(\xi_{i_1}))\cdots 
(\wh{\xi}_{i_r}+\gamma^\k(\xi_{i_r}))x $$
the result follows. 
\end{proof}

It is now easy to verify the following properties of the
Harish-Chandra map for Weil algebras.

\begin{theorem}[Harish-Chandra projection for Weil algebras] \label{th:hcweil}
Suppose $\g$ is a quadratic Lie algebra, 
and $\g=\n_-\oplus \k\oplus \n_+$ a decomposition
into subalgebras (direct sum of vector spaces) where $\k$ is 
quadratic and $\n_\pm$ are $\k$-invariant and isotropic. 
\begin{enumerate}
\item 
The maps 
$\kappa_W:\,W\g\to W\k$ and 
$\kappa_\W:\,\W\g\to \W\k$ are $\kds$ homomorphisms.
\item The diagram
$$
\xymatrix{W\g\ar[r]_{\ca{Q}_\g}
\ar[d]_{\kappa_{W}}&\W\g\ar[d]^{\kappa_{\W}}\\ 
W\k\ar[r]_{\ca{Q}_\k}       &\W\k}
$$
commutes up to $\k$-chain homotopy. 
\item
The above diagram contains a sub-diagram, 
$$ \xymatrix{S\g\ar[r] \ar[d]_{\kappa_S}&
U\g \ar[d]^{\tau\circ \kappa_U}\\ S\k\ar[r]
&U\k}$$
where $S\g$ is identified with $(W\g)_{\g-\on{hor}}$, 
$U\g$ with $(\W\g)_{\g-\on{hor}}$, and similarly 
for $S\k$ and $U\k$. In this sub-diagram the upper and 
lower horizontal map are Duflo maps, the left vertical map 
is $\kappa_S$, and the right vertical map is 
$\tau\circ \kappa_U$. 
\end{enumerate}
\end{theorem}

\begin{proof} 
It is clear that the decomposition \eqref{eq:decomp2} is 
$\wt{\k}$-equivariant, and that each summand is preserved 
by the differential. That is, \eqref{eq:decomp2} is a direct sum of 
$\k$-differential subspaces, which proves (a). 
Part (b) is immediate from Theorem \ref{th:rigidity}. 
We have already shown that the symmetrization
maps for Weil algebras restrict to the Duflo maps, and it is clear
that the left vertical map in (c) is just $\kappa_S$. 
The fact that the vertical map is  $\tau\circ \kappa_U$
follows from Proposition \ref{prop:kappas}.
\end{proof}
Theorem \ref{th:hcweil}
implies the following generalization of the
Harish-Chandra homomorphism for enveloping algebras.
\begin{theorem}\label{th:hc1}
The following diagram commutes: 
$$ \xymatrix{(S\g)_{\g-\on{inv}}\ar[r] \ar[d]_{\kappa_S}&
(U\g)_{\g-\on{inv}} \ar[d]^{\tau\circ \kappa_U}\\ 
(S\k)_{\k-\on{inv}} \ar[r]
&(U\k)_{\k-\on{inv}} }$$
Here the horizontal maps are the Duflo isomorphisms for 
$\g$ and $\k$, respectively.
\end{theorem}

\begin{proof}
By part (b) of the above Theorem, the diagram obtained by 
passing to  $\k$-basic cohomology 
$$ 
\xymatrix{{H((W\g)_{\k-\on{basic}})} \ar[r]
                \ar[d]&{H((\W\g)_{\k-\on{basic}})}
                \ar[d]\\
                (S\k)_{\k-\on{inv}}\ar[r]&(U\k)_{\k-\on{inv}}}
$$
commutes. (Moreover, all maps in this diagram are algebra isomorphisms.) 
On the other hand, the maps from $\g$-basic cohomology 
to $\k$-basic cohomology gives a commutative 
diagram of algebra homomorphisms,  
$$ 
\xymatrix{  (S\g)_{\g-\on{inv}}\ar[r]\ar[d]&(U\g)_{\g-\on{inv}}\ar[d] 
\\  
H((W\g)_{\k-\on{basic})}\ar[r]&H((\W\g)_{\k-\on{basic}})}
$$
Placing these two diagrams on top of each other, it follows that the
diagram in Theorem \ref{th:hc1} commutes.  
\end{proof}

\begin{proposition}
Under the assumptions of Theorem \ref{th:hc1}, the image of 
the cubic Dirac operator $\ca{D}_\g$ under the Harish-Chandra 
projection is 
$$ \kappa_\W(\ca{D}_\g)=\ca{D}_\k,$$
the cubic Dirac operator $\ca{D}_\k$ for the subalgebra. 
\end{proposition}

\begin{proof}
Recall that $\ca{D}_\g = \ca{D}_\k + \ca{D}_{\g, \k}$,
where $\ca{D}_{\g, \k} \in (U\g \otimes \Cl(\p))_{\k-\on{inv}}$. 
The image of $\ca{D}_{\g, \k}$ under the Harish-Chandra
projection vanishes since it is $\k$-basic and odd.
Hence, $\kappa_\W(\ca{D}_\g) = \kappa_\W(\ca{D}_\k)= \ca{D}_\k$. 
\end{proof}

\begin{remark}
For semi-simple Lie algebras and $\k=\h$ a Cartan subalgebra, the
Harish-Chandra projection $\kappa_{\on{Cl}}$ 
for Clifford algebras was studied by
Kostant. In particular Kostant showed that the image of a primitive
generator of $\wedge\g\cong \on{Cl}(\g)$ is always {\em linear},
i.e. contained in $\h\subset \on{Cl}(\h)$. He made a beautiful
conjecture relating these projections to the adjoint representation of
the principal TDS; this conjecture was recently proved by Y. Bazlov.
It would be interesting to understand these results within
our framework.  
\end{remark}

\section{Rouvi{\`e}re isomorphism}\label{sec:rouviere}
In his 1986 paper \cite{ro:es}, F. Rouvi{\`e}re described
generalizations of Duflo's isomorphism to a certain class of symmetric
spaces $G/K$. In this Section, we will prove a Duflo-Rouvi{\`e}re
isomorphism for quadratic Lie algebras $\g$, with a scalar product
that is anti-invariant under a given involution of $\g$.
\subsection{Statement of the theorem}
Let $\epsilon:\,\g\to \g$ be an involutive automorphism of a Lie algebra
$\g$. Then $\g=\k\oplus \mf{p}$ where $\k$ is the subalgebra fixed by
$\epsilon$, and $\mf{p}$ is the $-1$ eigenspace of $\epsilon$.  We will refer
to $(\g,\k)$ as a symmetric pair.  For any Lie algebra homomorphism
$f:\,\k\to \F$, define a twisted inclusion of $\k$ in $U\k$ by
$$ \k^f=\{\xi+f(\xi)|\,\xi\in\k\}.$$
Using the embedding $U\k\hra U\g$, we may view $\k^f$ as a subspace of 
$U\g$. The space 
\begin{equation}\label{eq:quot}
(U\g/ U\g\cdot \k^f)_{\k-\on{inv}}
\end{equation}
inherits an algebra structure from the
enveloping algebra $U\g$: Indeed, $U\g\cdot \k^f$ is a two-sided
ideal in the subalgebra $\{z\in U\g|\,L_\xi z\in U\g\cdot \k^f\mbox{
for all }\xi\in\k\}$ of $U\g$, and \eqref{eq:quot} is 
the quotient algebra. The following was proved by Duflo, generalizing 
a result of  Lichnerowicz \cite{li:op2}: 
\begin{theorem}[Duflo \cite{du:ope}] 
Let $(\g,\k)$ be a symmetric pair, and let $f:\,\k\to\F$ be the 
character $f(\xi)=\hh \on{tr}_\k(\ad_\xi)$. Then the 
algebra \eqref{eq:quot} is commutative.
\end{theorem}

The geometric interpretation of the algebra \eqref{eq:quot} is as
follows. Suppose $\F=\R$, and let $G$ be the connected, simply
connected Lie group having $\g$ as its Lie algebra. Assume that
$\k\subset\g$ is the Lie algebra of a closed, connected subgroup of
$K\subset G$. Taking $f=0$, a theorem of Lichnerowicz \cite{li:op1} shows
that \eqref{eq:quot} is the algebra of $G$-invariant differential
operators on the symmetric space $G/K$.  The algebra \eqref{eq:quot}
for $f(\xi)=\hh \on{tr}_\k(\ad_\xi)$ is interpreted as the algebra of
$G$-invariant differential operators on $G/K$, acting on sections of
the {\em half density bundle}.

Returning to the general case, we 
relate the algebra \eqref{eq:quot} to invariants in the symmetric 
algebra $S\mf{p}$. Indeed, using a PBW basis one sees that the map  
$$ S\mf{p}\oplus U\g\cdot \k^f\to U\g,\ \ (x,z)\mapsto \on{sym}_{U\g}(x)+z
$$
is a $\k$-module isomorphism. We therefore obtain an isomorphism  
of $\k$-modules
$$ \on{Sym}:\,S\mf{p}\to U\g/U\g\ \k^f$$
taking $x\in S\mf{p}$ to the image of $\on{sym}_{U\g}(x)$ under the
quotient map. Let $J_\mf{p}\in \ol{S\mf{p}^*}$ be defined by the function
$$ J_{\mf{p}}(\zeta)=\det(j(2\ad_\zeta)|_\mf{p})$$
with $j(z)=\f{\sinh(z/2)}{z/2}$ as in Section \ref{sec:duflo}.  This
is well-defined: For $\zeta\in\mf{p}$, $\ad_\zeta$ takes $\k$ to
$\mf{p}$ and vice versa; since $j$ is an even function, it follows
that $j(2\ad_\zeta)$ preserves both $\k$ and $\mf{p}$. 
Let
$$ \wh{J_{\mf{p}}^{1/2}}:\,S\mf{p}\to S\mf{p}$$
denote the infinite order differential operator defined by the
square root of the function $J_{\mf{p}}$. In Section 
\ref{subsec:proof1}  we will show: 
\begin{theorem}\label{th:canon}
Let $(\g,\k)$ be a symmetric pair, where $\k$ is the fixed point set 
of an involutive automorphism $\epsilon$. Suppose $\g$ admits an 
invariant scalar product $B$ with $\epsilon^*B=-B$. 
Then the composition
$$\on{Sym}\circ \wh{J_{\mf{p}}^{1/2}}:\, (S\mf{p})_{\k-\on{inv}}
\to (U\g/U\g\ \k^f)_{\k-\on{inv}}$$
(where $\mf{p}$ is the $-1$ eigenspace for $\epsilon$) is an algebra
isomorphism.
\end{theorem}
\begin{remark}
\begin{enumerate}
\item
A result similar to Theorem \ref{th:canon} was first proved by
Rouvi\`{e}re \cite{ro:es} for symmetric pairs $(\g,\k)$ satisfying one
of the following two conditions: (i) $\g$ is solvable, or (ii) $\g$
satisfies the Kashiwara-Vergne conjecture \cite{ka:ca} and $(\g,\k)$
is {\em very symmetric} in the sense that there is a linear
isomorphism $A: \g \to \g$, $A(\k)=\p$, $A(\p)=\k$ such that $[A,
\ad_\xi]=0$ for all $\xi\in \g$. At the time of this writing the
Kashiwara-Vergne conjecture is still open, but it has been established
for solvable Lie algebras and for quadratic Lie algebras \cite{ve:ce}
(see also \cite{al:ka}). (A major consequence of the
conjecture, regarding convolution of invariant distributions, 
was proved in the series of papers 
\cite{an:co,an:de,an:ko}. See \cite{to:ka} for more
information on the Kashiwara-Vergne method.) 
\item
It is known that for general symmetric pairs,
the statement of Theorem \ref{th:canon} becomes false. A
counter-example, examined in \cite{du:sim}, is $\g=\mf{sl}(2,\R)$,
with $\k=\so(1)$ the subalgebra of diagonal matrices, and $\mf{p}$ the
subspace of matrices having $0$ on the diagonal.
\end{enumerate}
\end{remark}

\subsection{Examples}\label{sec:examples}
Rouvi\`{e}res example $K^\C/K$ is included in the setting of Theorem
\ref{th:canon}, as follows.  Suppose $(\k,B_\k)$ is a quadratic Lie
algebra over $\F=\R$, and $\g=\k^\C$ is its complexification. Using
the extension of $B_\k$ to a complex-bilinear form on $\g$, define a
real-bilinear form
$$B_\g(\xi,\eta)=2\on{Im}(B_\k^\C(\xi,\eta))$$ 
on $\g$, viewed as a real Lie algebra. Then $B$ changes sign under the 
involution $\epsilon$ of $\g$ given by complex conjugation. 
The bilinear form $B_\k$ identifies $\k\cong\k^*$, while $B_\g$
identifies $\k^*\cong \mf{p}=\sqrt{-1}\k$.  The resulting map $\k\to
\sqrt{-1} \k$ is given by $\xi\mapsto \hh \sqrt{-1} \xi$. The function
$J_{\mf{p}}(\xi)$ turns into the following function
$$ J_h(\xi)=\det(j_h(\ad_\xi)),\ \ j_h(z)=\f{\sin(z/2)}{z/2},$$
similar to the usual Duflo factor, but with a $\sin$-function instead
of a $\on{sinh}$-function.  See \cite{al:li} for a geometric 
interpretation of this example. 

Suppose $\g$ is a quadratic Lie algebra, with involution $\eps$
changing the sign of the scalar product $B$. Then $B$ vanishes on both
$\k$ (the $+1$ eigenspace of $\eps$) and $\mf{p}$ (the $-1$ eigenspace
of $\eps$), and hence defines a non-singular pairing between $\k$ and
$\mf{p}$. This identifies $\mf{p}$ and $\k^*$, and defines an element
$C\in (\wedge^3\k)_{\k-\on{inv}}$ by 
$$
 B([\mu,\mu']_\g,\mu'')=C(\mu,\mu',\mu''),\ \ \mu,\,\mu',\,\mu''\in\k^*.
$$
Conversely, given a Lie algebra $\k$ and an invariant element 
$C\in \wedge^3\k$, the direct sum $\g=\k\oplus \k^*$ carries a 
unique Lie bracket such that $\k$ is a Lie subalgebra, 
$[\xi,\mu]_\g=-\ad_\xi^*\mu$ for $\xi\in \k,\,\mu\in\k^*$, and 
$$
 [\mu,\mu']_\g=C(\mu,\mu',\cdot)\in (\k^*)^*=\k,\ \ \mu,\,\mu'\in\k^*.
$$
Furthermore, the symmetric bilinear form $B$ given by the pairing
between $\k$ and $\k^*$ is $\g$-invariant, and changes sign under the
involution $\eps$ given by $-1$ on $\k^*$ and by $1$ on $\k$.
Hence all examples for Theorem \ref{th:canon} may be described in terms of 
a Lie algebra $\k$ with a given element $C\in (\wedge^3\k)_{\k-\on{inv}}$.

\begin{examples}
\begin{enumerate}
\item 
If $C=0$, the Lie algebra $\g$ is just the semi-direct product 
$\g=\k\ltimes\k^*$. In this case, one finds that 
$(U\g/U\g\ \k^f)_{\k-\on{inv}}=(S\k^*)_{\k-\on{inv}}$ and the 
Duflo-Rouvi\`{e}re isomorphism is just the identity map. 
\item
Suppose $(\k,B_\k)$ is a quadratic Lie algebra over $\F=\R$. 
Then $C(\xi,\xi',\xi'')=\pm B_\k([\xi,\xi']_\k,\xi'')$ defines an 
element $C\in (\wedge^3\k)_{\k-\on{inv}}$ (where we use $B_\k$ to identify 
$\k^*$ with $\k$). For the minus sign, one obtains the example 
$\g=\k^\C$ considered above. For the plus sign, one arrives at 
$\g=\k\oplus \k$, with $\k$ embedded diagonally and $\mf{p}$ embedded 
anti-diagonally. In this case, the Duflo-Rouvi\`{e}re isomorphism 
reduces to the usual Duflo isomorphism. 
\item 
Given $n\ge 3$ 
let $\k$ be the nilpotent Lie algebra of strictly upper-triangular $n\times
n$-matrices, and $\{E_{ij},\,i<j\}$ its natural basis, where $E_{ij}$
is the matrix having $1$ in the $(i,j)$ position and zeroes
elsewhere. Then $C=E_{1,n-1}\wedge E_{1,n}\wedge E_{2,n}$ lies in
$(\wedge^3\k)_{\k-\on{inv}}$. The resulting quadratic Lie algebra 
$\g$ is solvable for $n=3$, and nilpotent for $n\ge 4$. 
\item
There are non-trivial examples of $C\in (\wedge^3\k)_{\k-\on{inv}}$
such that the resulting symmetric pair $(\g,\k)$ is not very symmetric
in Rouvi\`{e}re's sense, and also $\g$ not solvable.  

Indeed if $(\g,\k)$ is a very symmetric pair with 
$[\k,\k]_\k=\k$, then  $[\mf{p},\mf{p}]_\g=\k$, or equivalently 
$\on{ker}(C)\equiv \{\mu\in\k^*|\,\iota_\mu C=0\}=0$. 
Take $\k=\a \ltimes \a^*$ with $\mf{a}$ semi-simple.  Let $C\in
\wedge^3\mf{a}^*\subset \wedge^3\k$ be defined by the Lie bracket and
the Killing form on $\mf{a}$. Since $C$ has non-zero kernel, $(\g,\k)$
is not very symmetric. Furthermore, $\a \subset \k \subset \g$ is the
Levi factor of $\g$ which shows that $\g$ is not solvable.
\end{enumerate}
\end{examples}

\subsection{Proof of Theorem \ref{th:canon}}
\label{subsec:proof1}
View $\W\g$ as a $\k$-differential algebra, with connection defined 
by the canonical $\g$-connection and the splitting $\g=\k\oplus\mf{p}$. 

Recall that the contraction operators on $\wt{\g}\subset \W\g$ are
given by $\iota_\xi\zeta=B(\xi,\zeta),\
\iota_\xi\ol{\zeta}=L_\xi\zeta$.  Since $\k$ is isotropic, it follows
that $\wt{\k}$ is a $\k$-differential subspace of $\W\g$, and so is
the left ideal $\W\g \wt{\k}$.  The algebra structure on $\W\g$ does
not descend to the $\kds$ $\W\g/(\W\g\ \wt{\k})$, in general. However,
there is an induced algebra structure on the basic subcomplex
$(\W\g/(\W\g\ \wt{\k}))_{\k-\on{basic}}$ since its pre-image in $\W\g$
is a subalgebra containing $\W\g\ \wt{\k}$ as a 2-sided ideal.


\begin{proposition}\label{prop:quot}
The characteristic homomorphism $W\k\to \W\g\cong U\g\otimes\Cl(\g)$ 
descends to a $\kds$ isomorphism 
\begin{equation}\label{eq:kdsisom}
W\k\cong S(E_\mf{p})\to \W\g/(\W\g\ \wt{\k})\cong 
(U\g/U\g\ \k^f)\otimes\wedge\mf{p},\end{equation}
an isomorphism of $\k$-modules 
\begin{equation}\label{eq:kmodisom}
 (W\k)_{\k-\on{hor}}\cong S\mf{p}\to (\W\g/(\W\g\
\wt{\k}))_{\k-\on{hor}}\cong U\g/U\g\ \k^f,
\end{equation}
and an isomorphism of algebras 
\begin{equation}\label{eq:algebraisom}
 (W\k)_{\k-\on{basic}}\cong (S\mf{p})_{\k-\on{inv}}
\to (\W\g/(\W\g\
\wt{\k}))_{\k-\on{basic}}\cong (U\g/U\g\ \k^f)_{\k-\on{inv}}.
\end{equation}
\end{proposition}
\begin{proof}
By a PBW argument, the map 
$$ S(E_\mf{p})\oplus \W\g\,\wt{\k}\to \W\g,\ (x,z)\mapsto Q_\g(x)+z$$
is a $\kds$ isomorphism. Thus, the quotient map $W\k\cong S(E_{\mf{p}})\to
\W\g/\W\g \wt{\k}$ is again a $\kds$ isomorphism, and its restriction to
horizontal subspaces is a $\k$-module isomorphism. 

The map on basic subspaces is a composition 
$$ (W\k)_{\k-\on{basic}}\to (\W\g)_{\k-\on{basic}}
\to (\W\g/(\W\g\ \wt{\k}))_{\k-\on{basic}},$$
where the second map is an algebra homomorphism, and the first map induces 
an algebra homomorphism in cohomology. Since the differential 
vanishes on $(W\k)_{\k-\on{basic}}$ 
and as a consequence vanishes on 
$(\W\g/(\W\g\ \wt{\k}))_{\k-\on{basic}}$, 
it follows that \eqref{eq:algebraisom} is an algebra isomorphism.
  
It remains to identify $\W\g/(\W\g\ \wt{\k})$ and its
horizontal and basic subspaces in terms of the isomorphism 
$\W\g=U\g\otimes \Cl(\g)$. 
Observe that $\W\g\ \wt{\k}$ is the
left ideal generated by elements $\zeta,\ol{\zeta}$ with
$\zeta\in\k$. We will show that
$$ \ol{\zeta}=\wh{\zeta}+\gamma^\g(\zeta)
=\wh{\zeta}+\hh \on{tr}_{\k}(\zeta) \mod \W\g\ \k.
$$
To see this, choose bases $e_i$ of $\k$ and $e^j$ of $\mf{p}$ such 
that $B(e_i,e^j)=\delta_i^j$. Then 
\beq \gamma^\g(\zeta)&=&\hh \sum_i (\ad_\zeta(e_i) e^i +\ad_\zeta(e^i) e_i)\\
&=&\hh \sum_i (-e^i\ad_\zeta(e_i) +\ad_\zeta(e^i) e_i)
+\hh \sum_i B(\ad_\zeta(e_i),e^i ).\eeq
The first sum lies in $\W\g\ \k$, while the second sum gives
$\hh \on{tr}_{\k}(\ad_\zeta)\mod \W\g\ \k$. 
This proves $\W\g/(\W\g\ \wt{\k})=(U\g/U\g\ \k^f)\otimes\wedge\mf{p}$,
where the contractions $\iota_\xi$ are induced by the contractions on
$\wedge\mf{p}=\wedge\k^*$. Hence the $\k$-horizontal subspace is
$(U\g/U\g\ \k^f)$, and the $\k$-basic subcomplex is $(U\g/U\g\
\k^f)_{\k-\on{inv}}$. 
\end{proof}
To complete the proof of Theorem \ref{th:canon}, we have to identify
the isomorphism \eqref{eq:kmodisom} from
$S\mf{p}$ onto $U\g/U\g \k^f$ with the map
$\on{Sym}\circ \wh{J_{\mf{p}}^{1/2}}$. Our calculation will require
the following Lemma.

\begin{lemma}\label{lem:contr}
Let $V$ be a vector space, and suppose $A:\,V\to V^*$ and $B:\,V^*\to V$ 
are linear maps with $A^*=-A,B^*=-B$. Let $\lambda(A)\in\wedge^2 V^*$ and 
$\lambda(B)\in\wedge^2 V$ be the skew-symmetric bilinear forms 
defined by $A,B$, i.e. in a basis $e_a$ of $V$, with dual basis 
$e^a$ of $V^*$, 
$$\lambda(A)=\hh \sum_a A(e_a)\wedge e^a,\ 
\lambda(B)=\hh \sum_a B(e^a)\wedge e_a.
$$
Let $\iota:\,V^*\to \on{End}(\wedge V)$ denote contraction, 
and denote by the same letter its extension to the exterior algebra 
$\wedge V^*$. Suppose $I+AB$ is invertible. Then 
$$ \iota(\exp(\lambda(A))\,\exp(\lambda(B))
={\det}^{1/2}(1+AB)\ \exp(\lambda(B\circ (I+AB)^{-1}))$$
for a unique choice of square root of $\det(I+AB)$. 
\end{lemma}
In particular, it follows that the map 
$$ \wedge^2 V^*\times \wedge^2 V\to \F,\ 
(\lambda(A),\lambda(B))\mapsto \det(I+AB)$$
admits a smooth square root, given by the degree zero part of
$\iota(\exp(\lambda(A))\,\exp(\lambda(B))$.
\begin{proof}
This may be proved by methods similar to \cite[Section 5]{al:cli}, to
which we refer for more details.  Let $V\oplus V^*$ be equipped with
the symmetric bilinear form given by the pairing between $V$ and
$V^*$, and $\on{Spin}(V\oplus V^*) \to \SO(V\oplus V^*)$ be the
corresponding Spin group. One has the following factorization in
$\SO(V\oplus V^*)$,
$$ \left(\begin{array}{cc}I&0\\A&I\end{array}\right)
\left(\begin{array}{cc}I&B\\0&I\end{array}\right)=
\left(\begin{array}{cc}I&D\\0&I\end{array}\right)
\left(\begin{array}{cc}I&0\\E&I\end{array}\right)
\left(\begin{array}{cc}R&0\\0&(R^{-1})^*\end{array}\right)
$$
where 
$$ R=(I+BA)^{-1},\ E=AR^{-1},\ D=BR^*.$$ 
This factorization lifts to a factorization in $\on{Spin}(V\oplus
V^*)$. Consider now the spinor representation of $\on{Spin}(V\oplus
V^*)$ on $\wedge V^*$. In this representation, the lift of  
{\tiny $ \Big(\!\!\begin{array}{cc}R&0\\0&(R^{-1})^*\end{array}\!\!\Big)$}
acts as $\alpha\mapsto {\det}^{-1/2}(R)\,R.\alpha$, the lift of 
{\tiny $ \Big(\!\!\begin{array}{cc}I&0\\E&I\end{array}\!\!\Big)$} acts
by contraction with $\exp(\lambda(E))$, and the lift of the factor 
{\tiny $ \Big(\!\!\begin{array}{cc}I&D\\0&I\end{array} \!\!\Big)$}
acts by exterior product with $\exp(\lambda(D))$. 
The lemma follows by applying the factorization to the ``vacuum vector''
$1\in\wedge V^*$. 
\end{proof}

\begin{proposition}\label{prop:whatitis}
The $\k$-module isomorphism $S\mf{p}
\to U\g/U\g\ \k^f$ in \eqref{eq:kmodisom}
is equal to the map 
$$\on{Sym}\circ \wh{J_{\mf{p}}^{1/2}}:\, S\mf{p}
\to U\g/U\g\ \k^f.$$
\end{proposition}
\begin{proof}
We write \eqref{eq:kmodisom} as a composition of maps 
$$
 S\mf{\p}
\stackrel{(i)}{\lra} W\g=S\g^*\otimes \wedge\g^*
\stackrel{(ii)}{\lra}\W\g=U\g\otimes\on{Cl}(\g) 
\stackrel{(iii)}{\lra} U\g
\stackrel{(iv)}{\lra} U\g/U\g\, \k^f
$$
Here (i) is the restriction of the characteristic map
$W\k\to W\g$ to $S\mf{\p}=S\k^*$, (ii) is the quantization map, 
(iii) is the tensor product of 
the augmentation map for $\on{Cl}(\g)$ with the identity map for $U\g$, 
and (iv) is the quotient map. 
In terms of the generators
$\wh{\mu}=\mu-\lambda^\k(\mu)$ of $S\k^*$, and the corresponding
generators of $S\g^*$, (i) is given by
$$ \mu\mapsto \mu,\ \ \wh{\mu}\mapsto \wh{\mu}+\lambda^{\mf{p}}(\mu)$$
where $\lambda^{\mf{p}}(\mu)=\lambda^\g(\mu)-\lambda^\k(\mu)$ takes values 
in $\wedge^2\mf{p}^*$. View $\lambda^{\mf{p}}$ as a 
$\wedge^2\mf{p}^*$-valued 
function on $\g^*$, constant in $\mf{p}^*$-directions. 
Then $\exp(\lambda^{\mf{p}})$ defines an infinite order 
$\wedge\mf{p}^*$-valued differential operator 
$$ \wh{\exp(\lambda^{\mf{p}})}:\,W\g\to W\g,$$
and the characteristic homomorphism $W\k\to W\g$ is a composition 
$\wh{\exp(\lambda^{\mf{p}})}\circ i$ where 
$i:\,W\k=S\k^*\otimes\wedge\k^*
\to W\g=S\g^*\otimes \wedge\g^*$ is the inclusion given on generators 
by $\mu\mapsto \mu,\wh{\mu}\mapsto \wh{\mu}$. 
Note that the image of $S\k^*\subset W\k$ under the composition 
lies in the 
subalgebra $S\k^*\otimes\wedge\mf{p}^*=S\mf{p}\otimes\wedge\k$ of 
$W\g$. Hence, when we apply the map (ii)
$$ \ca{Q}_\g=(\on{sym}\otimes q)\circ \iota(\wh{\S}):\,W\g\to \W\g,
$$
we need only consider the ``restriction'' of 
$\S(\xi)=J^{1/2}(\xi)\exp(\mf{r}(\xi)),\ \xi\in\g\cong\g^*$
to $\k^*\cong\mf{p}$. That is, we have to compute
\begin{equation}\label{eq:havetofind}
 J^{1/2}(\xi)\iota(\exp(\mf{r}(\xi)))\exp(\lambda^{\mf{p}}(\xi))
\in W\g=S\g\otimes\wedge\g
\end{equation}
for $\xi\in\mf{p}=\k^*$. In fact, we are only interested in the component 
of \eqref{eq:havetofind} in $S\g\otimes \wedge^0\g$, since all the 
other components will vanish under the projection (iii). Since 
$\k$ and $\mf{p}$ are dual $\k$-modules, 
$$ J(\xi)={\det}_{\mf{p}}(j(\ad_\xi)){\det}_\k(j(\ad_\xi))
={\det}_{\mf{p}}(j(\ad_\xi)^2).$$
Similarly, $r(\xi)$ splits into a sum $r'(\xi)+r''(\xi)$ where 
$r'(\xi)\in \wedge^2\mf{p}$ and $r''(\xi)\in \wedge^2\k$. 
We may replace $r(\xi)$ with $r'(\xi)$ in 
\eqref{eq:havetofind} since the $r''(\xi)$ part will not contribute 
to the contraction. Thus, we may calculate \eqref{eq:havetofind} 
using Lemma \ref{lem:contr}, with 
$$V=\mf{p},\ A=(\ln j)'(\ad_\xi)|_{\mf{p}},\ B=\ad_\xi|_{\mf{p}^*}.$$
Letting $(\cdot)_{[0]}$ denote the $\wedge^0\g$ component, 
the Lemma gives, 
$$
J^{1/2}(\xi)\,\big(\iota(\exp(\mf{r}(\xi)))\exp(\lambda^{\mf{p}}(\xi))\big)_{[0]}
={\det}^{1/2}_{\mf{p}}(h(\ad_\xi))$$
where 
\beq h(z)&=&j(z)^2 (1+z (\ln j)'(z))\\
&=& \big(\f{\on{sinh}(z/2)}{z/2}\big)^2
\,\big(1+z(\hh \on{coth}(z/2)-z^{-1})\big)\\
&=& \f{\on{sinh}(z/2)\on{cosh}(z/2)}{z/2}\\
&=&\f{\on{sinh}(z)}{z}.
\eeq
To summarize, the composition of the maps (i),(ii),(iii) is 
the map $S\mf{p}\to U\g$ given by application of an infinite-order 
differential operator $\wh{J_{\mf{p}}^{1/2}}$ on $S\mf{p}$, 
followed by the PBW symmetrization map $\on{sym}:\, S\mf{p}\to U\g$ . 
Since the map $\on{Sym}:\,S\mf{p}\to 
U\g/U\g\, \k^f$ is defined as PBW symmetrization
 followed by the quotient map (iv), the proof is complete. 
\end{proof}

\subsection{Isotropic subalgebras}\label{subsec:isotropic}
Theorem \ref{th:canon} may be generalized, as follows.  Suppose $\g$
is a quadratic Lie algebra, and $\k\subset \g$ an isotropic subalgebra
with $\k$-invariant complement $\mf{p}$.  Suppose $\mf{p}$ is {\em
co-isotropic}, that is, the $B$-orthogonal subspace $\mf{p}^\perp$
is contained in $\mf{p}$. 
Then $B$ induces a non-singular pairing between $\k$ and
$\mf{p}^\perp$, and the restriction of $B$ to 
$$\mf{m}=\k^\perp\cap \mf{p}$$ 
is non-singular. Since $\k$ is isotropic, $\wt{\k}\subset
\wt{\g}\oplus \F\mf{c}$ is a $\k$-differential subspace, and we may form
the quotient $\kds$, $\W\g/\W\g\, \wt{\k}$.  The proof of Proposition
\ref{prop:quot} carries over with no change, and shows that the map
$S(E_\mf{p})\to \W\g/\W\g\, \wt{\k}$ given by the inclusion
$S(E_\mf{p})\hra W\g$, followed by the quantization map and the
quotient map, is a $\kds$ isomorphism.

Since the scalar product identifies $\p^\perp \cong \k^*$, 
the $\kda$ $S(E_\mf{p})$ is of Weil type, with connection 
$\k^*\to E_\mf{p}$ given by the inclusion as 
$\mf{p}^\perp\subset \mf{p}$. It follows that the map in basic cohomology 
$$ S(\p^\perp)_{\k-\on{inv}}\to H((\W\g/\W\g\, \wt{\k})_{\k-\on{basic}})$$
is an algebra isomorphism.  The algebra on the right hand side admits
the following alternative description.  Let $\gamma^{\mf{m}}:\,\k\to
\on{Cl}\mf{m}$ be the map defined by the $\k$-action on $\mf{m}$. Then
$$ \gamma^\g(\zeta)=\hh \on{tr}_\k(\ad_\zeta)+\gamma^{\mf{m}}(\zeta)
\ \mod \W\g \wt{\k}.$$
Define a twisted inclusion $\k\to U\g\otimes \on{Cl}(\mf{m})$ 
by $ \zeta\mapsto \wh{\zeta}+ \hh \on{tr}_\k(\ad_\zeta)+
\gamma^{\mf{m}}(\zeta)$ and let $\k^f$ denote the image of this inclusion. 
As in the proof of Proposition \ref{prop:quot}, one finds that 
the algebra $(\W\g/\W\g\, \wt{\k})_{\k-\on{basic}}$ is canonically 
isomorphic to  the algebra 
\begin{equation}\label{eq:new}
\Big(\f{U\g\otimes
\on{Cl}(\mf{m})}{(U\g\otimes\on{Cl}(\mf{m}))\k^f}\Big)_{\k-\on{inv}}.
\end{equation}
One therefore arrives at the following result: 
\begin{theorem}\label{th:iso}
Let $\g$ be a quadratic Lie algebra, with isotropic subalgebra $\k$
and $\k$-invariant co-isotropic complement $\mf{p}$. Let 
$\mf{m}=\k^\perp\cap \mf{p}$. The algebra  \eqref{eq:new}
carries a natural differential, with cohomology algebra canonically 
isomorphic to $(S\k^*)_{\k-\on{inv}}$.  
\end{theorem}
Suppose $\F=\R$, and that $G$ is a connected Lie group
and $K$ a closed connected subgroup having $\g,\k$ as their Lie
algebras. Then \eqref{eq:new} has the following geometric 
interpretation. For any $K$-module $V$, consider the algebra $
\on{DO}(G\times_K V)_{G-\on{inv}}$ of $G$-invariant differential
operators on the associated bundle $G\times_K V$.  There is a
canonical algebra isomorphism (cf. \cite{li:op1})
$$ \on{DO}(G\times_K V)_{G-\on{inv}}=
\Big(\f{U\g\otimes \on{End}(V)}{(U\g\otimes
\on{End}(V))\,\k}\Big)_{K-\on{inv}}$$
where $\k\to U\g\otimes \on{End}(V)$ is embedded diagonally.  Thus, if
$\S$ is any module for the Clifford algebra $\Cl(\mf{m})$, and if the
action of the Lie algebra $\k$ (via $\k\to \Cl(\mf{m})$) on $\S$
exponentiates to an action of the Lie group $K$, there is an algebra
homomorphism from the algebra \eqref{eq:new} to the algebra of invariant 
differential operators on $G\times_K\ca{S}$, twisted by the half density bundle.
If the module $\ca{S}$ is faithful and irreducible, $\on{Cl}(\mf{m})\cong 
\on{End}(\ca{S})$ and this algebra homomorphism is an isomorphism.

\section{Universal characteristic forms}\label{sec:universal}
As a final application of our theory, we obtain a new construction of
universal characteristic forms in the Bott-Shulman complex
\cite{bo:le,du:si,mo:no,shu:th}.  We assume $\F=\R$, and let $G$ be a
connected Lie group with Lie algebra $\g$. Recall that in the
simplicial construction of the classifying bundle $EG\to BG$
\cite{mo:no,se:cl}, one models the de Rham complex of differential
forms on $EG$ by a double complex $C^{p,q}=\Om^q(G^{p+1})$ with
differentials $\d:\,C^{p,q}\to C^{p,q+1}$ (the de Rham differential)
and $\delta:\,C^{p,q}\to C^{p+1,q}$. Here $\delta$ is the alternating
sum $\delta=\sum (-1)^i\,\delta_i$ where $\delta_i=\partial_i^*$ is
the pull-back under the map,
$$\partial_i:\,G^{p+1}\to G^p,\ \ (g_0,\ldots,g_p)\mapsto 
(g_0,\ldots,\wh{g_i},\ldots,g_p)$$
omitting the $i$th entry. 
View each $E_pG=G^{p+1}$ as a principal $G$-bundle over
$B_pG=G^p$, with action the diagonal $G$-action from the right, and
quotient map $E_pG\to B_pG,\,(g_0,\ldots,g_p)\mapsto
(g_0g_1^{-1},g_1g_2^{-1},\ldots)$. Let 
$$ \iota_\xi:\,C^{p,q}\to C^{p,q-1},\ L_\xi:\,C^{p,q}\to C^{p,q}$$
be the corresponding contraction operators and Lie derivatives. 
Then the total complex  
$$ W=\bigoplus_{k=0}^\infty W^k,\ \ W^k=\bigoplus_{p+q=k} C^{p,q}$$
with differential $D=\d+(-1)^q\delta$, contractions $\iota_\xi$ and
Lie derivatives $L_\xi$ becomes a $\gds$. (The $\Z_2$-grading is given
by $W^\0=\bigoplus_{k=0}^\infty W^{2k}$ and
$W^\1=\bigoplus_{k=0}^\infty W^{2k+1}$.) By the simplicial de Rham
theorem \cite[Theorem 4.3]{mo:no}, 
the total cohomology of the basic subcomplex computes the
cohomology of the classifying space $BG$, with coefficients in $\R$.
Define a product structure on the double complex, 
$$ C^{p,q}\otimes C^{p',q'}\to C^{p+p',q}\otimes C^{p+p',q'}\to C^{p+p',q+q'}$$
where the first map is given by $(-1)^{qp'}$ times the tensor products
of pull-backs to the first $p+1$, respectively last $p'+1$, $G$-factors in
$G^{p+p'+1}$, and the second map is wedge product. (This formula 
is motivated by the usual formula for cup products of the singular 
cochain complex.) It is
straightforward to verify that $D,\iota_\xi,L_\xi$ are derivations for
the product structure, thus $W$ becomes a $\gda$. It is locally free,
with a natural connection
$$\theta:\,\g^*\to C^{0,1}=\Om^1(G)$$ 
given by the left-invariant Maurer-Cartan form on $G$. Hence, by 
symmetrization we 
obtain a $\gds$ homomorphism $ W\g\to  W$ which
restricts to a map of basic subcomplexes. The resulting map
\begin{equation}\label{eq:bottshulmann}
 (S\g^*)_{\on{inv}}\to \bigoplus_{p,q}\Om^q(G^{p+1})_{\on{basic}}
\cong \bigoplus_{p,q} \Om^q(G^p)
\end{equation}
takes an invariant polynomial of degree $r$ to a $D$-cocycle of total
degree $2r$. By our general theory, the induced map in cohomology is a
ring homomorphism. As in the usual Bott-Shulman construction we have 
the following vanishing phenomenon:
\begin{proposition}
The image of an invariant polynomial of degree $r$ under the map 
\eqref{eq:bottshulmann} has non-vanishing components only in bidegrees
$(p,q)$ with $p+q=2r$ and $p\le r$. 
\end{proposition}
\begin{proof}
The connection $\theta$ lives in bidegree $(0,1)$, and its 
total differential $D\theta$ has non-vanishing components 
only in bi-degrees $(0,2)$ and $(1,1)$. Since the product structure 
is compatible with the bi-grading, it follows that the image of 
an element $\ol{\xi}_{i_1}\cdots \ol{\xi}_{i_l}\xi_{j_1}\cdots 
\xi_{j_m}\in W\g$ under the symmetrization map only 
involves bi-degrees $(p,q)$ with $p+q=2l+m$ and $p\le l$. 
Any element in $S^r\g^*\subset W\g$ is a linear combination of such elements, 
with $2l+m=2r$.  
\end{proof}

Let $W^L\subset W$ denote the direct sum over the subspaces 
$(C^{p,q})^L=\Om^q(G^{p+1})^L$ of forms that are 
invariant under the left $G^{p+1}$-action. Clearly, $W^L$ is 
a $\g$-differential subalgebra of $W$. Since the connection 
is left-invariant, our Chern-Weil map takes values in the 
subalgebra $W^L_{\on{basic}}$.  
\begin{theorem}
$W^L$ is a $\gda$ of Weil type. If $G$ is compact and connected, the
$\gda$ $W$ is of Weil type, with the inclusion $W^L\hra W$ a
$\g$-homotopy equivalence.
\end{theorem}
It follows that in this case, the map \eqref{eq:bottshulmann} gives an
isomorphism in cohomology.
\begin{proof}
We first recall the standard proof that $W$ is acyclic.  Let
$\Pi:\,W\to W$ be the projection operator given on $C^{p,q}$ as
pull-back under the map
$$ \pi:\,G^{p+1}\to G^{p+1},\ (g_0,\ldots,g_p)\mapsto (e,\ldots,e).$$
 The projection $\Pi$ is naturally chain homotopic to the identity:
To construct a homotopy operator, let
$$ s_j:\,G^{p+1}\to  G^{p+2},\ \ 
(g_0,\ldots,g_{p})\mapsto 
(g_0,\ldots,g_{j},e,\ldots,e),\ \ i=0,\ldots,p$$ 
and set $s=\sum_{j=0}^p (-1)^j s_j^*:\,C^{p+1,q}\to C^{p,q}$.  Note
that on $G^{p+1}$, $ \partial_{p+1}s_p=\on{id},\ \partial_0 s_0=\pi$.
A direct calculation (as in May \cite{ma:si}) shows that
$[\delta,s]=\Pi-\on{id}$ on $\bigoplus_{p=0}^\infty\Om^q(G^{p+1})$,
for any fixed $q$. Since $[\d,s]=0$, it follows that $-(-1)^q
s:\,\Om^q(G^{p+1})\to \Om^q(G^{p+2})$ gives the desired homotopy
between the identity and $\Pi$. 

The image $\Pi(W)\subset W$ is isomorphic to the singular cochain
complex of a point. Hence, composing with the standard  
homotopy operator for this
complex we see that the inclusion $i:\,\F\hra W$ is a homotopy 
equivalence. Let $h:\,W\to W$ denote the homotopy operator 
constructed in this way. 

The map $s$ does not commute with $L_\xi$ since the maps $s_j$ are not
$G$-equivariant. However, on the left invariant subcomplex $W^L$ this
problem disappears, since $s_j^*\circ L_\xi=L_\xi\circ s_j^*$ on
$W^L$. It follows that $h$ restricts to a homotopy operator on $W^L$
with $[h,L_\xi]=0$.  This shows that $W^L$ is of Weil type.

Suppose now that $G$ is compact. Then there is a projection
$\Pi_1:\,W\to W^L$, given on $\Om^q(G^{p+1})$ as the averaging
operator for the left $G^{p+1}$-action. It is well-known that the
averaging operator is homotopic to the identity operator; the homotopy
operator $h_1$ may be chosen to commute with $L_\xi$, by averaging
under the right $G$-action. (In our case, one may directly construct
$h_1$ using the Hodge decomposition for the bi-invariant Riemannian
metric on $G$.)  It follows that $\Pi|_{W^L}\circ \Pi_1$ is homotopic
to the identity map, by a homotopy operator that commutes with $L_\xi$
and lowers the total degree by $1$.
\end{proof}
\begin{remark}
For the classical groups, there is another model for differential
forms on the classifying bundle, as an inverse limit of differential
forms on Stiefel manifolds -- e.g. for $G=\U(k)$, the inverse limit of
$\Om(\on{St}(k,n))$ for $n\to \infty$.  The resulting $\gda$ carries a
natural ``universal'' connection (see Narasimhan-Ramanan
\cite{na:ex}).  The characteristic homomorphism for this case was
studied by Kumar in \cite{ku:re}.
\end{remark}

\begin{appendix}
\section{Proof of Theorem \ref{th:technical}}
Recall that a linear operator $C$ on a vector space $E$ is {\em
locally nilpotent} if $E=\bigcup_{N\ge 0} E^{(N)}$ where $E^{(N)}$ is
the kernel of $C^{N+1}$. If $C$ is locally nilpotent, the operator
$I+C$ has a well-defined inverse, since the geometric series
$I-C+C^2-C^3\pm \cdots$ is finite on each $E^{(N)}$.
\begin{lemma}\label{lem:boot}
Suppose $E$ is a $\g$-differential space, and $h\in\on{End}(E)^\1$ 
an odd linear operator with the following properties: 
$$ [\iota_\xi,h]=0,\ [L_\xi,h]=0,\ [\d,h]=I-\Pi+C$$
where $\Pi$ is a projection operator, and $C$ is locally nilpotent. 
Assume that $h$ and $C$ vanish 
on  the range of $\Pi$. Then $\ti{\Pi}=\Pi(I+C)^{-1}$ is a
projection operator having the same range as $\Pi$. 
Furthermore, it is 
a $\gds$ homomorphism, and 
$\ti{h}=h (I+C)^{-1}$ is a $\g$-homotopy between $I$ and $\ti\Pi$:
$$ [\iota_\xi,\ti{h}]=0,\ [L_\xi,\ti{h}]=0,\ [\d,\ti{h}]=I-\ti{\Pi}.$$
\end{lemma}
\begin{proof}
Since $C\Pi=0$, it is clear that $\ti{\Pi}$ is a projection
operator. We will check $[\d,\ti{h}]=I-\ti{\Pi}$ on any 
given $v\in E$. Choose $N$ sufficiently large so that 
$v\in E^{(N)}$ and $\d v\in E^{(N)}$. 
On $E^{(N)}$, the operator $\ti{h}$ is a finite series
\beq \ti{h}&=&h(I-C+C^2- \cdots+(-1)^N C^N)
\\&=&h(I-(C-\Pi)+(C-\Pi)^2+ \cdots
+(-1)^N\,(C-\Pi)^N)
,\eeq
where we have used $h\Pi=0$ and $C\Pi=0$. 
The equation $[\d,h]=I-\Pi+C$ shows that $C-\Pi$ commutes with 
$\d$. Thus
$$ [\d,\ti{h}]=[\d,h](I-(C-\Pi)+(C-\Pi)^2+ \cdots
+(-1)^N\,(C-\Pi)^N)$$
on the subspace $\{v\in E^{(N)}|\,\d v\in E^{(N)}\}$. 
But $[\d,h]=I-\Pi+C$ also shows $[\d,h]\Pi=0$. Hence we may replace
$C-\Pi$ by $C$ again, and get 
\beq [\d,\ti{h}]\,&=&[\d,h](I-C+C^2+\cdots)\\
&=&[\d,h](I+C)^{-1}\\
&=&(I-\Pi+C)(I+C)^{-1}=I-\ti{\Pi}. 
\eeq
By a similar argument, since $[L_\xi,C-\Pi]=[L_\xi,[h,\d]]=0$ 
and 
$[\iota_\xi,C-\Pi]=-[\iota_\xi,[h,\d]]
=-[h,L_\xi]=0,$
one proves $[L_\xi,\ti{h}]=0$ and $[\iota_\xi,\ti{h}]=0$. 
\end{proof}

\begin{lemma}\label{lem:conn2}
Suppose $E$ is a $\g$-differential space, and $E'\subset E_{\on{basic}}$ 
a differential subspace. Suppose $h:\,E\to E$ is an odd linear
operator with $[L_\xi,h]=0$ and $[\d,h]=I-\Pi$, where $\Pi$ is a
projection operator onto $E'$ with $h\Pi=0$. 
Assume there exists an increasing
filtration $E=\bigcup_{N=0}^\infty E^{(N)}$, with $E'\subset
E^{(0)}$, such that each $\iota_\xi$ and $h$ have 
negative filtration degree. Then, for any
locally free $\g$-differential space $\ca{B}$, the inclusion map
$$E'\otimes\ca{B}\to E\otimes\ca{B}$$
is a $\g$-homotopy equivalence, with a homotopy inverse
$E \otimes\ca{B}\to E'\otimes\ca{B}$ that is equal to the identity
on $E'\otimes\ca{B}$.  
\end{lemma}
\begin{proof}
The proof is inspired by an argument of Guillemin-Sternberg (see
\cite[Theorem 4.3.1]{gu:su}).  If $[\iota_\xi^E,h]=0$, the projection
operator $\Pi=\Pi\otimes I$ is the desired homotopy inverse to the
inclusion map, with homotopy $h=h\otimes 1$. In the general case, we
employ the Kalkman trick \cite{ka:br} to shift the contraction
operators on $E\otimes\ca{B}$ to the second factor.  Let $e_a$ be a
basis of $\g$, and $e^a$ the dual basis of $\g^*$.  Choose a
connection $\theta:\,\g^*\to \ca{B}$.  Then $\psi=\sum_a
\theta(e^a)\iota_{e_a}^E$ is an even, nilpotent operator on
$E\otimes\ca{B}$, with $[\psi,L_\xi^E+L_\xi^{\ca{B}}]=0$.  Hence
$\exp \psi$ is a well defined even, invertible operator on
$E\otimes\ca{B}$, and it commutes with the $\g$-action. Let
$$ 
\ti{L}_\xi=\Ad(\exp\psi)(L^E_\xi+L^\ca{B}_\xi),\ \ 
\ti{\iota}_\xi=\Ad(\exp\psi)(\iota^E_\xi+\iota^\ca{B}_\xi),\ \ 
\ti{\d}=\Ad(\exp\psi)(\d^E+\d^\ca{B})$$
denote the transformed Lie derivatives, contractions an differential
on $E\otimes\ca{B}$.  A calculation using
$\Ad(\exp\psi)=\exp(\ad\psi)=\sum_{n=0}^\infty \f{1}{n!}(\ad\psi)^n$
shows
\beq
\ti{\iota}_\xi&=&\iota^{\ca{B}}_\xi\\
\ti{L}_\xi&=&L^E_\xi+L^{\ca{B}}_\xi\\
\ti{\d}&=&\d^E+\d^{\ca{B}}
+\theta(e^a)L_{e_a}^E+R \eeq
where the remainder term $R$ is a polynomial in contractions
$\iota_{e_a}^E$, with coefficients in $\ca{B}$ (with no constant 
term). The operator $h$ commutes with the 
contractions $\ti{\iota}_\xi$ and Lie derivatives $\ti{L}_\xi$, 
and 
$$ [\ti{\d},h]=[\d^E,h]+[R,h]=I-\Pi+[R,h].$$
Our assumptions imply that $C=[R,h]$ has negative filtration degree,
and in particular is locally nilpotent.  Thus Lemma \ref{lem:boot}
applies and shows that $\ti{\Pi}=\Pi(I+C)^{-1}$ is a projection
operator with the same range $E'\otimes \ca{B}$, and is homotopic to
the identity by a homotopy $\ti{h}=h(I+C)^{-1}$ which commutes with
all $\ti{\iota}_\xi,\ \ti{L}_\xi$.  Since $\iota_\xi^E$ vanishes on
$E'$, the operator $\exp\psi$ acts trivially on $E'\otimes\ca{B}$ and
therefore $\wh{\Pi}=\ti{\Pi}\circ \exp(\psi)$ is a projection onto
$E'\otimes\ca{B}$.  The operator $\wh{h}=\Ad(\exp(-\psi))\ti{h}$ gives
the desired $\g$-homotopy between $\wh{\Pi}$ and the identity.
\end{proof}

\begin{proof}[Proof of Theorem \ref{th:technical}]
Let $W,W'$ be $\gda$'s of Weil type. Lemma \ref{lem:conn2} (for
$E=W'$ and $\ca{B}=W$) gives a $\g$-homotopy equivalence $W\to
W'\otimes W,\ \ w\mapsto 1\otimes w$.  Reversing the roles of $W,W'$,
the Lemma also shows that there exists a $\g$-homotopy equivalence
$W\otimes W'\to W'$ that restricts to the identity map on $\F\otimes
W'$.  By composition, we obtain a $\g$-homotopy equivalence $W\to W'$
taking units to units.
\end{proof}

\end{appendix}

\bibliographystyle{amsplain}   %


\def\polhk#1{\setbox0=\hbox{#1}{\ooalign{\hidewidth
  \lower1.5ex\hbox{`}\hidewidth\crcr\unhbox0}}} \def\cprime{$'$}
  \def\cprime{$'$} \def\polhk#1{\setbox0=\hbox{#1}{\ooalign{\hidewidth
  \lower1.5ex\hbox{`}\hidewidth\crcr\unhbox0}}} \def\cprime{$'$}
\providecommand{\bysame}{\leavevmode\hbox to3em{\hrulefill}\thinspace}
\providecommand{\MR}{\relax\ifhmode\unskip\space\fi MR }
\providecommand{\MRhref}[2]{%
  \href{http://www.ams.org/mathscinet-getitem?mr=#1}{#2}
}
\providecommand{\href}[2]{#2}

\end{document}